\documentclass[a4paper,11pt]{article}

\usepackage{amsmath,amsthm,amssymb,amsfonts}
\usepackage{mathtools} 
\usepackage{bm}
\usepackage{mathrsfs} 
\usepackage{changes} 
\usepackage{url} 
\usepackage[toc,page]{appendix} 
\usepackage{enumerate}

\usepackage{graphicx} 
\usepackage{longtable} 
\usepackage[all]{xy} 
\usepackage{newfloat} 
\DeclareFloatingEnvironment[fileext=frm,placement={ht},name=Frame]{algfloat}
\usepackage{caption}
\usepackage{placeins}
\usepackage{tikz}
\usepackage{subfig}
\usepackage{color}
\usepackage{xcolor} 
\usepackage{array}
\usepackage{multirow}
\usepackage{booktabs}
\usepackage{grffile}

\newcommand{\N}{\mathbb{N}}

\newcommand{\R}{\mathbb{R}}

\newcommand{\inner}[2]{\langle #1,#2\rangle} 

\newcommand{\norm}[1]{\|{#1}\|}
\newcommand{\normq}[1]{\|{#1}\|^2}

\newcommand{\HH}{{\mathcal{H}}} 

\newcommand{\tos}{\rightrightarrows} 

\newcommand{\comenta}[1]{} 

\newcommand{\set}[1]{\{#1\}}

\newcommand{\lab}[1]{\label{#1}}

\newcommand{\bprop}{\begin{proposition}}
\newcommand{\eprop}{\end{proposition}}
\newcommand{\blemm}{\begin{lemma}}
\newcommand{\elemm}{\end{lemma}}
\newcommand{\bdefi}{\begin{definition}}
\newcommand{\edefi}{\end{definition}}
\newcommand{\btheo}{\begin{theorem}}
\newcommand{\etheo}{\end{theorem}}
\newcommand{\bproo}{\begin{proof}}
\newcommand{\eproo}{\end{proof}}
\newcommand{\brema}{\begin{remark}}
\newcommand{\erema}{\end{remark}}
\newcommand{\bitem}{\begin{itemize}}
\newcommand{\eitem}{\end{itemize}}
\newcommand{\bexam}{\begin{example}}
\newcommand{\eexam}{\end{example}}
\newcommand{\bassu}{\begin{assumption}}
\newcommand{\eassu}{\end{assumption}}
\newcommand{\bcoro}{\begin{corollary}}
\newcommand{\ecoro}{\end{corollary}}

\newcommand{\mgap}{\vspace{.1in}}

\newcommand{\beq}{\begin{equation}}
\newcommand{\eeq}{\end{equation}}
\newtheorem{theorem}{Theorem}[section] 
\newtheorem{lemma}[theorem]{Lemma}
\newtheorem{corollary}[theorem]{Corollary}
\newtheorem{proposition}[theorem]{Proposition}
\newtheorem{remark}[theorem]{Remark}
\newtheorem{definition}[theorem]{Definition}
\newtheorem{assumption}[theorem]{Assumption}
\newtheorem{example}[theorem]{Example}

\usepackage[linesnumbered,ruled]{algorithm2e}

\SetKwInput{KwInput}{Input}
\SetKwInput{KwOutput}{Global Variables for Function}
\SetKw{KwGlob}{Global Variables for Function:}

\SetKwFunction{FMain}{BrackBisecProc} 
\SetKwProg{Fn}{Function}{:}{}

\usepackage[T1]{fontenc}
\usepackage[utf8]{inputenc}
\usepackage[english]{babel} 

\definechangesauthor[name={M. Marques Alves}, color=blue]{M} 
\definechangesauthor[name={B. Svaiter}, color=purple]{B}

\begin{document}

\title{A search-free $O(1/k^{3/2})$ homotopy inexact proximal-Newton extragradient algorithm for monotone variational inequalities} 
\author{
    M. Marques Alves%
\thanks{
Departamento de Matem\'atica,
Universidade Federal de Santa Catarina,
Florian\'opolis, Brazil, 88040-900 ({\tt maicon.alves@ufsc.br}).
Partially supported by CNPq grant 308036/2021-2.
}
\and
J. M. Pereira%
\thanks{ IMPA, Estrada Dona Castorina 110, 22460-320 Rio de
    Janeiro, Brazil ({\tt jpereira@impa.br}). }
\and
B. F. Svaiter%
\thanks{ IMPA, Estrada Dona Castorina 110, 22460-320 Rio de
    Janeiro, Brazil ({\tt benar@impa.br}). 
    Partially supported by CNPq grants 430868/2018-9 and 311300/2020-0.}
}


\maketitle

\begin{abstract}
We present and study the iteration-complexity of a relative-error inexact
proximal-Newton extragradient algorithm for solving smooth monotone variational
inequality problems in real Hilbert spaces. 
We removed a search procedure from Monteiro and Svaiter (2012) by introducing a novel approach based on
homotopy, which requires the resolution (at each iteration) of a single strongly monotone
linear variational inequality.
For a given tolerance $\rho>0$, our main algorithm exhibits 
pointwise $O\left(\frac{1}{\rho}\right)$ and ergodic
$O\left(\frac{1}{\rho^{2/3}}\right)$ iteration-complexities.
From a practical perspective, preliminary numerical experiments indicate that our main algorithm outperforms some previous proximal-Newton schemes.
\end{abstract}

\pagestyle{plain}


\section{Introduction}\lab{sec:intr}

This paper is concerned with the numerical solution of the (monotone) variational inequality problem (VIP)
\begin{align} \lab{eq:probi}
\mbox{Find}\;\;x\in C\;\;\mbox{such that}\;\; \inner{F(x)}{y-x}\geq 0\qquad \forall y\in C,
\end{align}
where $F(\cdot)$ is a smooth monotone operator with Lipschitz continuous derivative $F'(\cdot)$
and $C$ is a closed and convex set (more details in Section \ref{sec:main}).
For examples and applications of \eqref{eq:probi} in different contexts in optimization and applied mathematics we refer the
reader, e.g., to \cite{fac.pan-fin01, fac.pan-fin02, roc.wet-book}.

We propose and study the iteration-complexity of an inexact proximal-Newton extragradient algorithm (see Algorithm \ref{alg:main} below) for the numerical solution of \eqref{eq:probi}. 
For the purposes of this paper, it will be convenient to reformulate problem \eqref{eq:probi} as a  monotone inclusion problem
\begin{align} \lab{eq:probim}
0\in F(x) + N_C(x),
\end{align}
where $N_C$ denotes the normal cone operator of $C$ (in the sense of convex analysis). 
It is a matter of fact that
\eqref{eq:probi} and \eqref{eq:probim} are equivalent, and so from now on in this paper we will refer to \eqref{eq:probim} as our main problem.

One iteration of a proximal Newton-type method for solving \eqref{eq:probim}, at a current step $x_{k-1}$, can be 
performed by solving
\begin{align}\lab{eq:reg.pro}
0\in \lambda_k \big(F_{x_{k-1}}(x) + N_C(x)\big) + x - x_{k-1},
\end{align}
where $\lambda_k>0$ and $F_{x_{k-1}}(\cdot)$ 
is a linearization of $F(\cdot)$ at $x_{k-1}$ as in \eqref{eq:def.lin} below (with $y=x_{k-1}$).
A proximal-Newton extragradient method for monotone inclusions and variational inequalities with iteration-complexity 
$O\left(\frac{1}{\rho^{2/3}} \log\log \frac{1}{\rho}\right)$ (where $\rho>0$ is a given precision) was proposed and analyzed in \cite{mon.sva-newton.siam12}.
Recently, the latter method was studied and generalized in different directions; see, e.g., 
\cite{bullins2022lai, jiang2022mokhtari, Lin2023Jordan}.
The \emph{multiplicative} $\log\log \frac{1}{\rho}$ factor which appears in the complexity analysis of the main algorithms in 
\cite{bullins2022lai, jiang2022mokhtari, Lin2023Jordan,mon.sva-newton.siam12} is the consequence of a search procedure
needed (at each iteration) in order to compute the solution of subproblems and for the tuning of the proximal parameters.
It has been posed as a challenging question how to remove 
the just mentioned search procedures (see, e.g., \cite[Section 4.3.3]{Nesterov2018Book}). 
At the cost of solving a nonlinear VIP
at each iteration, significant progress was made in this direction in \cite{adil.bullins.jambulapati2022, jordan-perseus22} (see \cite[Definition 3.1]{adil.bullins.jambulapati2022} and Eq. (7) in \cite{jordan-perseus22}).

In this paper, we propose a different approach to address this
issue. We propose an inexact proximal-Newton extragradient algorithm based on ideas from homotopy, without destroying the subproblem formulation originally proposed in \cite{mon.sva-newton.siam12}. More precisely, our main contributions can be summarized as follows.

\mgap

\noindent
{\bf Main contributions.} 
\begin{itemize}
\item[(i)] Our main algorithm is Algorithm \ref{alg:main} below, which is a relative-error (inexact) proximal-Newton extragradient method for solving \eqref{eq:probim}. We emphasize that, in contrast to \cite{adil.bullins.jambulapati2022, jordan-perseus22}, 
Algorithm \ref{alg:main} preserves the same subproblem structure originally proposed in 
\cite{mon.sva-newton.siam12}, which is to say that at each iteration it
requires the inexact solution of a strongly monotone \emph{linear} VIP defined
by the linearized operator -- see \eqref{eq:err.cri}.
Moreover, a single instance of this class of VIP has to be inexactly solved at each iteration.
See also the third comment following Algorithm \ref{alg:main} and 
Subsections \ref{subsec:subprob-a} and \ref{subsec:subprob} for more details on the solution of subproblems.
\item[(ii)] Algorithm \ref{alg:main} is a search-free method
with  pointwise 
$O\left(\frac{1}{\rho}\right)$ and ergodic
$O\left(\frac{1}{\rho^{2/3}}\right)$ iteration-complexities;
see Theorems \ref{th:main}, \ref{th:main02} and Corollary \ref{cr:bigo}. 
Summarizing, our main algorithm is a search-free inexact proximal-Newton extragradient method with 
pointwise $O\left(\frac{1}{\rho}\right)$ and ergodic $O\left(\frac{1}{\rho^{2/3}}\right)$ iteration-complexities, respectively, which preserves the (prox-regularized) subproblem structure from \cite{mon.sva-newton.siam12}.
\end{itemize}

\mgap

\noindent
{\bf Some related works.}
Papers~\cite{bullins2022lai, jiang2022mokhtari, Lin2023Jordan} generalized the proximal-Newton extragradient (NPE) method of Monteiro and Svaiter~\cite{mon.sva-newton.siam12} to more general settings. 
More precisely, \cite{bullins2022lai} proposed and studied a high-order mirror-prox method that achieves an ergodic $O\left(1/k^{(p+1)/2}\right)$ global convergence rate ($p=2$ gives the ergodic rate $O\left(1/k^{3/2}\right)$ of~\cite{mon.sva-newton.siam12}). 
The same ergodic rate for $p\geq 2$ was also obtained
in~\cite{jiang2022mokhtari} for a generalized optimistic method for convex-concave saddle-points, where the case of strongly-convex-strongly-concave problems was also considered.
Paper~\cite{Lin2023Jordan} proved pointwise $O\left(1/k^{p/2}\right)$ and ergodic $O\left(1/k^{(p+1)/2}\right)$ (global) convergence rates (in both continuous and discrete regimes) for a $p^{\emph{th}}$-order generalization of the NPE method, with
a local linear convergence rate also studied under an appropriate error-bound condition.
We mention that in the aforementioned works, namely \cite{bullins2022lai, jiang2022mokhtari, Lin2023Jordan}, all generalizations of the NPE method of Monteiro and Svaiter also depend on a search procedure for computing approximate solutions of subproblems. 
As we mentioned before, this search procedure was removed in~\cite{adil.bullins.jambulapati2022, jordan-perseus22} at the cost of solving more complex subproblems. 
In~\cite{adil.bullins.jambulapati2022}, an ergodic $O\left(1/k^{(p+1)/2}\right)$ rate was proved for a variant of the mirror-prox method of Nemirovski in noneuclidean settings. 
For solving variational inequalities on bounded domains, \cite{jordan-perseus22} obtained 
the same ergodic  $O\left(1/k^{(p+1)/2}\right)$ rate for a variant of the dual extrapolation method of Nesterov as well as global linear and local superlinear rates under strong monotonicity assumptions. 
Both works~\cite{adil.bullins.jambulapati2022, jordan-perseus22} proved that their methods are optimal by computing the lower bound
$\Omega(\varepsilon^{-2/(p+1)})$. 
Motivated by the (continuous-time) perspective of acceleration as the discretization of differential equations, 
\cite{Lin2023AContinuous} proposed and studied novel accelerated rescaled gradient systems for solving nonlinear monotone equations in the form of $F(x) = 0$. Moreover, first- and high-order methods achieving a global rate of 
$O(k^{-p/2})$ in terms of the residual $\norm{F(x)}$ are obtained from the discretization of the continuous-time system.
For nonconvex and nonconcave structured minmax optimization, \cite{Vyas2023Beyond} presented a high-order variant of the extragradient method, which can be derived as a discrete time version of the rescaled gradient 
system of \cite{Lin2023AContinuous} (see Algorithm 1 in \cite{Vyas2023Beyond}). For the main algorithm 
of \cite{Vyas2023Beyond}, the authors showed that it finds $\varepsilon$-appoximate stationary  at a rate of $O(k^{-p/2})$; the same rate was also obtained in the continous-time setting.
Paper \cite{Lin2023Explicit} obtained an ergodic $O\left(1/k^{3/2}\right)$ iteration-complexity for an explicit 
second-order algorithm for convex-concave unconstrained saddle-point problems; their main algorithm is a search-free cubic-regularized Newton method \cite[Algorithm 1]{Lin2023Explicit}. In the latter reference, inexact variants (for Hessian evaluation and subproblems solving) of the main algorithm are also studied with applications to finite-sum minmax problems with subsampling of the Hessian. At this point we emphasize again that, although search-free, contrary to our work at each iteration \cite[Algorithm 1]{Lin2023Explicit} demands the resolution of more complex cubic-regularized subproblems. 
On the other hand, \cite{Lin2023Explicit} contains a characterization of the actual final runtime of subproblem solving.

Complexity-based second- and high-order methods for convex optimization have also been the subject of intense research in 
the last two decades (see, e.g., \cite{Arjevani2019Oracle, Bubeck2019Near-optimal, Gasnikov2019Optimal, Jiang2021Anoptimal, Monteiro2013Accelerated, Nesterov2021Implementable, Nesterov2008Accelerating, Nesterov2006Cubic}). 
Recently, the concurrent papers~\cite{Carmon2022Optimal, Kovalev2022TheFirst} proposed and studied search-free 
optimal second- and high-order methods for convex optimization.  
In~\cite{Carmon2022Optimal}, optimal and adaptive search-free variants of the Monteiro-Svaiter (MS) method are proposed; a subproblem solver for the solution of linear systems based on exact solution or MinRes is also studied and numerical experiments on logistic regression problems are presented. 
On the other hand, paper~\cite{Kovalev2022TheFirst} proposed an optimal (search-free) tensor extragradient method for smooth unconstrained convex optimization.
We emphasize that although search-free, the algorithms 
discussed in the aforementioned references are tailored especifically for convex optimization and are based on different 
ideas and goals, when compared to our work.

\mgap

\noindent
{\bf Our approach.} 
In the context of this paper, the main step in the NPE method of \cite{mon.sva-newton.siam12} consists in computing an inexact solution $y_k = y_k(\lambda_k, x_{k-1})$ of \eqref{eq:reg.pro} satisfying
\begin{align} \lab{eq:tea}
 \dfrac{2\sigma_\ell}{L}\leq \lambda_k\norm{y_k - x_{k-1}}\leq \dfrac{2\sigma_u}{L},
\end{align}
where $\sigma_u > \sigma_\ell >0$ and $L > 0$ is the Lipschitz constant of $F'(\cdot)$. Note that, since $y_k$ depends on both
$\lambda_k > 0$ and $x_{k-1}$, the proximal parameter $\lambda_k>0$ and the approximate solution $y_k$ have to be properly tuned in order to satisfy the (coupled) system \eqref{eq:reg.pro}--\eqref{eq:tea}. 
This problem was solved for the first time in \cite{mon.sva-newton.siam12} by presenting and analyzing the complexity of a search-procedure for computing $(\lambda_k, y_k)\in \R_+\times\HH$ satisfying \eqref{eq:reg.pro}--\eqref{eq:tea} ($y_k$ satisfying \eqref{eq:reg.pro} as an inexact solution). 
As we discussed before, this search-procedure has been replicated and explored in various research papers. A key question that has emerged is how to eliminate this procedure, and this challenge serves as one of the primary motivations for this work.

With this goal in mind, in this paper we explore some techniques from 
homotopy (or continuation methods), which have been used in interior-point 
methods~\cite{Nesterov2018Book, Nesterov1994Book}, in the setting of variational inequalities and proximal-type algorithms.
To give some context of the main ideas behind our main algorithm, we first note that the two inequalites in \eqref{eq:tea} have distinct goals. The second inequality is related to the progress of the iterations produced by 
the NPE method \cite{mon.sva-newton.siam12} toward the solution set of \eqref{eq:probim}, whereas the first one, which we now write as (cf. condition $\lambda_k\norm{y_k - x_{k-1}}\geq \eta$ in Algorithm \ref{alg:main})
\begin{align}\lab{eq:book}
 \lambda_k\norm{y_k - x_{k-1}}\geq \dfrac{2\sigma_\ell}{L},
\end{align} 
is to enforce larger steps along the iterative process, enabling the optimal convergence rate $O(1/k^{3/2})$ 
(see \cite[Theorem 2.5]{mon.sva-newton.siam12}). In our approach, the right-hand side of \eqref{eq:book} will be replaced by a fixed parameter $\eta > 0$, which depends on the Lipschitz constant $L > 0$ in a different way (see Eq. \eqref{eq:def.tau.alg} below).
Then our main algorithm (Algorithm \ref{alg:main} below) proceeds by  computing an approximate solution of \eqref{eq:reg.pro} but now with $F(\cdot)$ linearized at a different point $y_{k-1}$ (see Eq. \eqref{eq:err.cri}). Next the algorithm checks if the large-step condition $\lambda_k\norm{y_k - x_{k-1}}\geq \eta$ holds. 
If this is the case, in order to make progress toward the solution set, an under-relaxed extragradient step is 
performed and the proximal parameter $\lambda_k>0$ is updated so that $\lambda_{k+1}<\lambda_k$ ; otherwise, we increase the proximal parameter to $\lambda_{k+1}>\lambda_k$ (see lns. 6 -- 10 of Algorithm \ref{alg:main}). For the steps such that the large-step condition $\lambda_k\norm{y_k - x_{k-1}}\geq \eta$ holds, the complexity analysis follows from the corresponding one for the large-step under-relaxed HPE method (see Section \ref{sec:basic} and Theorem \ref{th:rhpe}); on the other hand, when the large-step condition fails, a separate analysis is needed (see Proposition \ref{prop:iteb} below).

\mgap

\noindent
{\bf General notation.}
Let $\HH$ be a real Hilbert space with inner product $\inner{\cdot}{\cdot}$ and induced norm 
$\norm{\cdot} = \sqrt{\inner{\cdot}{\cdot}}$.
We use the same notation $\norm{A}$ to denote the norm of a bounded linear operator $A$.
The \emph{normal cone} to $C\subset \HH$ at $x\in C$ is defined as
%
$N_C(x) = \left\{\nu\in \HH\;\;|\;\;\inner{\nu}{y-x}\leq 0\quad \forall y\in C\right\}.$
%
We also set $N_C(x) := \emptyset$ if $x\notin C$. The $\varepsilon$-enlargement of the set-valued map $T:\HH\tos \HH$ (see, e.g., \cite{bur.ius.sva-enl.svva97}) is the (set-valued map) $T^\varepsilon:\HH\tos \HH$ defined as
\begin{align}\lab{eq:def.enl}
T^{\varepsilon}(x) = \left\{v\in \HH\;\;|\;\;\inner{v - u}{x - z}\geq -\varepsilon\quad \forall u\in T(z)\right\}. 
\end{align}
We also use the notation 
\begin{align}\lab{eq:def.lin}
F_y(x) := F(y) + F'(y)(x-y)\qquad \forall x\in C
\end{align}
to indicate the linearization of a smooth map $F:C\subset \HH\to \HH$ at $y\in C$.
The set of natural numbers is $\N = \{1,2,\dots\}$ and $\# S$ stands for the number of elements of
a set $S$. We also define $\log^+(t) = \max\{\log(t), 0\}$ for $t>0$.

\mgap

\section{On the large-step under-relaxed HPE method} \lab{sec:basic}

In this section, we present a special instance of the large-step 
under-relaxed hybrid proximal-extragradient (HPE) method~\cite{pre-print-benar} (see Algorithm \ref{alg_gen} below) and its iteration-complexity analysis; see Theorem \ref{th:rhpe}. These results will be instrumental in the complexity analysis of Algorithm \ref{alg:main} (see Proposition \ref{prop:special}).

\mgap
\mgap

\begin{algorithm}[H]
\lab{alg:rHPE}
 \caption{A special instance of the large-step under-relaxed HPE method for \eqref{eq:probim}}
\lab{alg_gen}
\SetAlgoLined
\KwInput{$x_0\in \HH$, $0<\tau<1$, $0\leq \sigma <1$, $\eta>0$ and $N\in \N\cup \{\infty\}$}
 \For{$k = 1,2,\dots, N$}{
  Choose $\lambda_k>0$ and find $y_k\in \HH$ and $\nu_k\in N_C(y_k)$ such that\\
  \begin{align}
  \lab{eq:rel.al1}
   &\left\|\lambda_k \left(F(y_k) + \nu_k\right) + y_k - x_{k-1}\right\|\leq \sigma\norm{y_k - x_{k-1}}\\[2mm]
   \lab{eq:lar.ste}
   &\lambda_k\norm{y_k - x_{k-1}}\geq \eta
  \end{align}
  Set
  \begin{align}\lab{eq:ext.al1}
   x_k = x_{k-1} - \tau\lambda_k \left(F(y_k) + \nu_k\right)
  \end{align}
  }
\end{algorithm}

\mgap
\mgap

We now make the following remarks regarding Algorithm \ref{alg:rHPE}:
\begin{itemize}
\item[(i)] Algorithm \ref{alg_gen} is a special instance of the large-step under-relaxed method of \cite{pre-print-benar} (which was proposed for general maximal monotone operators) for solving \eqref{eq:probim}. 
If we set $\tau = 1$ in \eqref{eq:ext.al1}, then it follows that Algorithm \ref{alg_gen} 
becomes a special instance of the large-step HPE method~\cite{mon.sva-newton.siam12} for solving \eqref{eq:probim}. 
We note that \eqref{eq:rel.al1} is a \emph{relative-error} condition for the inexact solution of the proximal subproblem
\begin{align}\lab{eq:exa.inc}
0\in \lambda_k\left( F(y) + N_C(y)\right) + y - x_{k-1}
\end{align}
in the sense that if we take $\sigma = 0$ in \eqref{eq:rel.al1}, then we have that $y =y_k$ is the 
solution of \eqref{eq:exa.inc}, i.e., $y_k = \left(\lambda_k(F+N_C) + I\right)^{-1}x_{k-1}$, and $\nu_k$ is given explicitly by
\[
 \nu_k = \dfrac{x_{k-1} - y_k}{\lambda_k} -  F(y_k).
\] 
On the other hand, depending on the particular structure of $F(\cdot)$ and $C$ in \eqref{eq:probim}, condition \eqref{eq:rel.al1} can be used as a stopping criterion for running an ``inner algorithm'' for computing approximate solutions
of \eqref{eq:exa.inc}. 
\item[(ii)] The large-step condition \eqref{eq:lar.ste} first appeared in \cite{mon.sva-newton.siam12}; it is tailored for the design and analysis of first-order methods for VIPs (or second-order methods for optimization) and, in particular, implies better complexity results (faster convergence rates) for the large-step HPE method 
when compared to the HPE method with constant stepsize (see~\cite{mon.sva-newton.siam12}). 
We mention that if we remove condition \eqref{eq:lar.ste} and set $\tau =1$ in Algorithm \ref{alg_gen}, then it reduces to the HPE method~\cite{mon.sva-hpe.siam10, sol.sva-hpe.svva99} for solving \eqref{eq:probim}.  
We also mention that condition \eqref{eq:ext.al1} is an under-relaxed extragradient step performed in order to update $x_k$.
\item[(iii)] Algorithm \ref{alg_gen} will be used as a framework for the design and complexity analysis of our main algorithm, namely Algorithm \ref{alg:main} below. The sequences of Algorithm \ref{alg:main} corresponding to indexes $k\geq 1$ satisfying a large-step condition can be regarded as an algorithm which turns out to be a special instance of Algorithm 
\ref{alg_gen} (see Proposition \ref{prop:special}).
\end{itemize}

\mgap

Before presenting the iteration-complexity of Algorithm \ref{alg_gen}, we need to clarify what do we mean by an approximate solution of \eqref{eq:probim}. For a given tolerance $\rho>0$, we want to find $y\in C$ such that either there exists $\nu\in N_C(y)$ such that
\begin{align}\lab{eq:app.solai}
\norm{F(y) + \nu}\leq \rho
\end{align}
or there exists $\varepsilon\geq 0$ and $v\in (F + N_C)^{\varepsilon}(y)$ such that
\begin{align}\lab{eq:app.solbi}
\max \left\{\norm{v},\varepsilon\right\}\leq \rho,
\end{align}
where $(F + N_C)^{\varepsilon}$ denotes the $\varepsilon$-enlargement of $F + N_C$ (see Section \ref{sec:intr} for general notation).
Conditions \eqref{eq:app.solai} and \eqref{eq:app.solbi} are closely related to the more usual notions of weak and strong solution of VIP problems (see, e.g., \cite{mon.sva-hpe.siam10} for a discussion).

\mgap

The \emph{ergodic sequences} associated to the sequences generated by Algorithm \ref{alg:rHPE} are defined as 
follows (see, e.g., \cite{mon.sva-hpe.siam10}):

\begin{align}\lab{eq:def.erg1}
& y_k^a = \dfrac{1}{\sum_{i=1}^k\,\lambda_i}\,\sum_{i=1}^k\,\lambda_i  y_i,\\[2mm]
\lab{eq:def.erg2}
& v_k^a = \dfrac{1}{\sum_{i=1}^k\,\lambda_i}\,\sum_{i=1}^k\,\lambda_i  \left(F(y_i) + \nu_i\right),\\[2mm]
\lab{eq:def.erg3}
& \varepsilon_k^a = \dfrac{1}{\sum_{i=1}^k\,\lambda_i}\,\sum_{i=1}^k\,\lambda_i \inner{y_i - y_k^a}{F(y_i) + \nu_i- v_k^a}.
\end{align}

\mgap

Next we present the convergence rates of Algorithm \ref{alg_gen}.

\mgap

\begin{theorem} 
 \lab{th:rhpe}
Consider the sequences evolved by \emph{Algorithm \ref{alg:rHPE}} and the ergodic sequences as in \eqref{eq:def.erg1}--\eqref{eq:def.erg3}. Let also $d_0$ denote the distance of $x_0$ to the solution set $(F + N_C)^{-1}(0)\neq \emptyset$ of \eqref{eq:probim}. The following hold for all $1\leq k\leq N$:
\begin{itemize}
\item[\emph{(a)}] There exists $1\leq j\leq k$ such that $\nu_j\in N_C(y_j)$ and
\begin{align*}
\norm{F(y_j) + \nu_j}\leq \dfrac{d_0^2}{\tau \eta (1-\sigma)\,k}\,.
\end{align*}
\item[\emph{(b)}] $v_k^a\in (F + N_C)^{\varepsilon_k^a}(y_k^a)$ and
\begin{align*}
\norm{v_k^a}&\leq \dfrac{2d_0^2}{\tau^{3/2}\eta\sqrt{1-\sigma^2}\,k^{3/2}}\,,\\[3mm]
\varepsilon_k^a&\leq \dfrac{2d_0^3}{\tau^{3/2}\eta(1-\sigma^2)\,k^{3/2}}.
\end{align*}
\end{itemize}
\end{theorem}
\bproo
The proof follows from \cite[Theorem 2.4]{pre-print-benar}) applied to the maximal monotone operator 
$T = F + N_C$.
\eproo

\mgap

\brema
\emph{
Theorem \ref{th:rhpe}(a) shows that, for a given tolerance $\rho>0$, one can find $y\in C$ and $\nu\in N_C(y)$ such that
\eqref{eq:app.solai} holds in at most
\begin{align*}
\left\lceil\left(\dfrac{1}{\tau\eta(1 - \sigma)}\right)\dfrac{d_0^2}{\rho}\right\rceil
\end{align*}
iterations. On the other hand, Theorem \ref{th:rhpe}(b) gives that one can find $y\in C$, $v\in \HH$ and $\varepsilon\geq 0$ such that
$v\in (F + N_C)^{\varepsilon}(y)$ and \eqref{eq:app.solbi} holds in at most 
\begin{align*}
\begin{aligned}
& \max\left\{\left\lceil{\left(\dfrac{\sqrt[3]{4}}
      {\tau \eta^{2/3}(1-\sigma^2)^{1/3}}\right)\left(\dfrac{d_0^2}{\rho}\right)^{2/3}}\right\rceil, 
      \left\lceil{\left(\dfrac{\sqrt[3]{4}}
      {\tau \eta^{2/3}(1-\sigma^2)^{2/3}}\right)\left(\dfrac{d_0^3}{\rho}\right)^{2/3}}\right\rceil\right\}
\end{aligned}
\end{align*}
iterations.
}
\erema

\mgap

\section{The main algorithm}\lab{sec:main}

Consider the (monotone) VIP given in \eqref{eq:probi}, i.e.,
\begin{align} \lab{eq:prob}
\mbox{Find}\;\;x\in C\;\;\mbox{such that}\;\; \inner{F(x)}{y-x}\geq 0\qquad \forall y\in C,
\end{align}
where $C$ is a (nonempty) closed and convex subset of $\HH$ and 
$F:C\to \HH$ is assumed to be (point-to-point) \emph{monotone}
and continuously differentiable with a $L$-Lipschitz continuous derivative, i.e., $F$ and $F'$ satisfy the following conditions for all
$x,y\in C$:
\begin{align}\lab{eq:mon.lip}
\begin{aligned}
& \inner{F(x) - F(y)}{x - y}\geq 0,\\[2mm]
& \norm{F'(x) - F'(y)}\leq L\norm{x-y}.
\end{aligned}
\end{align}

Recall that \eqref{eq:prob} is equivalent to \eqref{eq:probim}, i.e.,
\begin{align}\lab{eq:probm}
 0 \in F(x) + N_C(x).
\end{align}
We assume throughout this paper that the solution set $(F+N_C)^{-1}(0)$ of \eqref{eq:probm} 
(or, equivalently, of \eqref{eq:prob}) is nonempty.

As we mentioned before, our main goal in this paper is to propose and study the iteration-complexity of an inexact 
proximal-Newton algorithm for solving \eqref{eq:probm} under the assumptions described in
\eqref{eq:mon.lip}. 

\mgap
\mgap

\noindent
{\bf Subproblem approximate solution.} One iteration of the proximal-Newton algorithm for solving \eqref{eq:probm}, at a current step $x_{k-1}$, can be performed by solving the inclusion
\begin{align*}
0\in \lambda_k \big(F_{x_{k-1}}(x) + N_C(x)\big) + x - x_{k-1},
\end{align*}
where $\lambda_k>0$ and $F_{x_{k-1}}(\cdot)$ is as in \eqref{eq:def.lin} (with $y=x_{k-1}$).
Our main algorithm 
(Algorithm \ref{alg:main} below) will be defined by solving, at each iteration, a similar inclusion:
\begin{align}\lab{eq:farias}
0\in \lambda_k \big(F_{y_{k-1}}(x) + N_C(x)\big) + x - x_{k-1},
\end{align}
where $(y_k)$ is an auxiliary sequence. Note that when $C=\HH$, in which case \eqref{eq:probm} reduces to the nonlinear equation $F(x)=0$, then solving \eqref{eq:farias} reduces to solving the linear system
\begin{align}\lab{eq:linsys}
\left(\lambda_k F'(y_{k-1}) + I\right)(x - y_{k-1}) = -\left(\lambda_k F(y_{k-1}) + y_{k-1} - x_{k-1}\right).
\end{align}
We refer the reader to the several comments following Algorithm \ref{alg:main} for additional discussions and remarks concerning the resolution of \eqref{eq:farias}.

\mgap
\mgap

We will allow inexact solutions of \eqref{eq:farias} within relative-error tolerances, which we define next.

\mgap

\bdefi[$\hat\sigma$-approximate solution of \eqref{eq:farias}]
 \lab{def:app_sol}
For $y_{k-1}\in C$, $x_{k-1}\in \HH$ and $\lambda_k>0$, a pair $(y_k,\nu_k)\in \HH\times \HH$ is said to be a 
$\hat\sigma$-approximate solution of \eqref{eq:farias} 
if $\hat\sigma\in [0,1/2)$ and
\begin{align}\lab{eq:app_sol}
\begin{aligned}
\begin{cases}
\nu_k\in N_C(y_k),\\[3mm]
\norm{\lambda_k\left(F_{y_{k-1}}(y_k) + \nu_k\right) + y_k - x_{k-1}} \leq \hat\sigma\norm{y_k - y_{k-1}}.
\end{cases}
\end{aligned}
\end{align}
\edefi


Let 
\begin{align}
  \lab{eq:hsigma}
  0 \leq \hat\sigma <1/2
\end{align}
be the parameter of approximate solution of \eqref{eq:farias}, as
in Definition~\ref{def:app_sol}.
This $\hat\sigma$ can be regarded as the inexactness of an (inexact) oracle which
solves approximately the strongly monotone linear VIP~\eqref{eq:farias}.

Before stating the main algorithm, 
choose and define the parameters $\theta$, $\hat\theta$, $\eta$, $\tau$
as follows:
\begin{align}\lab{eq:def.hth}
  0 < \theta<(1-\hat\sigma)(1-2\hat\sigma),
  \quad 
  \hat\theta & := \theta\left(\dfrac{\hat\sigma}{1-\hat\sigma}
+ \dfrac{\theta}{(1-\hat\sigma)^2}\right)
\end{align}
and 
\begin{align}\lab{eq:def.tau.alg}
    \eta & > \dfrac{2\hat\theta}{L},
  \quad \tau  := \dfrac{2 \left(\theta - \hat\theta\right)}
{2\theta + \dfrac{\eta L}{2} + \sqrt{\left(2\theta + \dfrac{\eta L}{2}\right)^2 
- 4\theta\left(\theta - \hat\theta\right)}}.
\end{align}

\mgap
\mgap

Next is our main algorithm.

\mgap
\mgap

\begin{algorithm}[H]
 \caption{A relative-error inexact proximal-Newton algorithm for solving \eqref{eq:probm}}
\lab{alg:main}
\SetAlgoLined
\KwInput{$x_0\in C$,$\;y_0:=x_0$, $\nu_0 := 0$, the scalars 
 $\hat\sigma$, $\theta$, $\;\hat\theta$ $\eta$, $\tau$
 as in \eqref{eq:hsigma}--\eqref{eq:def.tau.alg}
and $\lambda_1>0$ such that $\lambda_1^2\norm{F(y_0)}\leq 2\theta/L$
}
 \For{$k=1,2,\dots$}{
 \If{$F(y_{k-1})+\nu_{k-1}=0$}
 	{\Return $y_{k-1}$\\\lab{eq:return}}
 \eIf{$\dfrac{\lambda_k L}{2}\norm{\lambda_k \left(F(y_{k-1}) + \nu_{k-1}\right) + y_{k-1} - x_{k-1}}\leq 
  \hat\theta$}
  {$y_k = y_{k-1}$\\ \lab{eq:n.cond}
   $\nu_k =\nu_{k-1}$ \\ \lab{eq:n.cond02}}
  {Find a $\hat\sigma$-approximate solution $(y_k,\nu_k)$ of \eqref{eq:farias} (in the sense of
   Definition \ref{def:app_sol}), i.e., find $(y_k,\nu_k)$ satisfying
   \begin{align}\lab{eq:err.cri}
    \begin{cases}
   \nu_k\in N_C(y_k),\\[3mm]
   \norm{\lambda_k\left(F_{y_{k-1}}(y_k) + \nu_k\right) + y_k - x_{k-1}} \leq \hat\sigma\norm{y_k - y_{k-1}}
   \end{cases}
  \end{align} \lab{eq:do.newton}}
   \eIf{$\lambda_k\norm{y_k - x_{k-1}}\geq \eta$}
      {
      \vspace{0.2cm}
      $x_{k} = x_{k-1} -\tau\lambda_k \left(F(y_k) + \nu_k\right)$\\[2mm]\lab{eq:ext.step}
      \vspace{0.1cm}
      $\lambda_{k+1} = (1 - \tau)\lambda_k$\\[3mm]\lab{eq:red.lam}
      }
     {$x_k = x_{k-1}$\\[2mm]\lab{eq:null.step}
     $\lambda_{k+1} = \dfrac{\lambda_k}{1 - \tau}$\\\lab{eq:inc.lam}
     }
  }
\end{algorithm}

\mgap
\mgap

Next we make some remarks regarding Algorithm \ref{alg:main}:

\begin{itemize}
\item[(i)] Note that since $y_0\in C$ and $\nu_0=0$, it follows that $\nu_0\in N_C(y_0)$. 
Note also that the condition on $\lambda_1>0$ is trivially satisfied if $F(y_0)=0$ and reduces to 
$0<\lambda_1\leq \sqrt{2\theta/(L\norm{F(y_0)})}$ otherwise. 
\item[(ii)] We mention that if Algorithm \ref{alg:main} stops at step \ref{eq:return}, then it follows that $y_{k-1}$ is a solution of
\eqref{eq:probm}. So from now on in this paper, we assume that
\[
\fbox{
Algorithm \ref{alg:main} doesn't stop at step \ref{eq:return}.
}
\]
\item[(iii)] As we mentioned before, if $C=\HH$ then \eqref{eq:farias} reduces to the linear system \eqref{eq:linsys}. In this case, the error criterion as given in \eqref{eq:err.cri} can be written in terms of the residual of \eqref{eq:linsys}:
\begin{align}\lab{eq:friday}
\norm{\left(\lambda_k F'(y_{k-1}) + I\right)(y_k - y_{k-1}) + \left(\lambda_k F(y_{k-1}) + y_{k-1} - x_{k-1}\right)}
\leq \hat\sigma\norm{y_{k} - y_{k-1}}.
\end{align}
So condition \eqref{eq:friday} can be used as a stopping rule for computing inexact solutions of \eqref{eq:linsys} by means
of an ``inner algorithm'' depending on the specific structure of $F$ (e.g., the conjugate gradient algorithm in the case of quadratic problems).
Moreover, if $F=\nabla f$, where $f:\HH\to \R$ is a twice continuously differentiable convex function (satisfying \eqref{eq:mon.lip} with $F'=\nabla^2f$), then problem \eqref{eq:probm} reduces to the
constrained minimization problem $\min_{x\in C}\,f(x)$ and \eqref{eq:farias} is clearly equivalent to solving the
regularized second-order model:
\[
\min_{x\in C}\,\left(f(y_{k-1}) + \inner{\nabla f(y_{k-1})}{x - y_{k-1}} + 
\dfrac{1}{2}\inner{\nabla^2f(y_{k-1})(x - y_{k-1})}{x - y_{k-1}} + \dfrac{1}{2\lambda_k}\normq{x - x_{k-1}}\right).
\]
We will also discuss the cost of solving the subproblems \eqref{eq:err.cri} in the general case in
Subsections \ref{subsec:subprob-a} and \ref{subsec:subprob} below.

\item[(iv)] The extragradient step as in line \eqref{eq:ext.step} is performed if the large-step condition $\lambda_k\norm{y_k - x_{k-1}}\geq \eta$ is satisfied; otherwise, we set $x_k = x_{k-1}$ and increase the value of $\lambda_k$ (see lines \ref{eq:null.step} and \ref{eq:inc.lam}).
\item[(v)] The iteration-complexity of Algorithm \ref{alg:main} will be presented in 
Theorems \ref{th:main} and \ref{th:main02} below (see also Corollary \ref{cor:rate} and Proposition \ref{prop:iteb}). The complexities are pointwise 
$\mathcal{O}\left(\frac{1}{\rho}\right)$
and ergodic $\mathcal{O}\left(\frac{1}{\rho^{2/3}}\right)$ 
(see, e.g., \eqref{eq:def.m} and \eqref{eq:def.mt}, respectively).
\end{itemize}

\mgap

\brema \lab{rem:theta}
\emph{We mention that
\begin{align}\lab{eq:tau.ineq}
  0<\hat\theta<\theta\;\;\;\mbox{and}\;\;\;0<\dfrac{\theta -\hat\theta}{2\theta + \dfrac{\eta L}{2}}<\tau<1.
\end{align}
The fact that $\theta>\hat \theta$ is a direct consequence of the assumption $0<\theta<(1-\hat\sigma)(1-2\hat\sigma)$
(see the input in Algorithm \ref{alg:main}) and \eqref{eq:def.hth}. On the other hand, by defining the quadratic function 
$q:\R\to \R$,
\begin{align}\lab{eq:def.q}
q(t) := \theta t^2 - \left(2\theta + \dfrac{\eta L}{2}\right)t + \theta -\hat\theta,\qquad t\in \R,
\end{align}
we see that $q(0) = \theta -\hat\theta >0$ and $q(1) = -\frac{\eta L}{2} - \hat\theta<0$. Then, since $\tau$ is clearly the smallest root of $q(\cdot)$, it follows that
\begin{align}\lab{eq:root.q}
 q(\tau) =0
\end{align}
and
\begin{align*}
-\dfrac{q(0)}{q'(0)}<\tau< 1,
\end{align*}
which is equivalent to \eqref{eq:tau.ineq}.}
\erema

\mgap

	\subsection{On the solution of subproblems} 
	\lab{subsec:subprob-a}
	At each iteration $k\geq 1$, Algorithm \ref{alg:main} -- see Eqs. \eqref{eq:farias} and \eqref{eq:err.cri} --, demands the
	computation of a $\hat\sigma$-approximate solution of the monotone inclusion
	\begin{align} \lab{eq:walt}
		0\in \lambda_k\left(F_{y_{k-1}}(y) + N_C(y)\right) + y - x_{k-1}.
	\end{align}
	Note that the above inclusion is equivalent to a variational inequality for the \emph{affine} \emph{strongly-monotone} operator
	$G_k: \HH \to \HH$,
	\begin{align} \lab{def:opg}
		G_k(y) := \lambda_k F_{y_{k-1}}(y) + y - x_{k-1},\qquad y\in \HH,
	\end{align}
	that is, to find (the unique) $y_k^*\in C$ such that
	\begin{align}
		\lab{eq:subp}
		\inner{G_k(y_k^*)}{y - y_k^*}\geq 0\qquad \mbox{for all}\;\;y\in C.
	\end{align}
	We propose several options for solving this problem, depending on
	properties and assumptions on $ C $:
	\\
	\\
	{\bf [a] } $ C = \HH $. In this case, the problem to be solved is a linear
	system. These problems may be solved efficiently and up to numerical precision using direct methods, such as LU decomposition, or inexactly using iterative methods, such as conjugate gradients, MinRes \cite{paige1975solution} or LSQR \cite{paige1982lsqr}.
	\\
	{\bf [b] } $ C $ is the positive orthant, the cone of symmetric positive semidefinite matrices,
	a second order cone or a convex set for which a self-concordant barrier is
	available.
	In this case, one may use an interior point method tailored to the problem at
	hand.
	%
	\\
	{\bf [c] } $ C = \set{x \mid g_i(x) \leq 0,\; \; i=1,\dots,m} $ where
	$ g_1,\dots,g_m $ are continuously differentiable and have Lipschitz
	continuous Hessians.  In this case, we propose a reformulation of the original problem:
	\begin{align*}
		\widetilde C = \HH \times \R_+^m,
		\qquad
		\widetilde F (x,\xi)= 
		\begin{bmatrix}
			F(x)+\sum_{i=1}^m \xi_i \nabla g_i(x)
			\\
			-g(x)
		\end{bmatrix}.
	\end{align*}
	Then, one may use an interior point method to solve the linearized proximal subproblems of 
	$ VIP (\widetilde F, \widetilde C ) $.
	\\
	{\bf [d]} We have only an oracle which computes the orthogonal projection onto $ C $.
	This setting is analyzed in the next subsection.

\mgap

	\subsection{On the solution of subproblems for the projection oracle model} 
	\lab{subsec:subprob}
	We assume in this subsection that $ C $ is accessed by means of a \emph{projection oracle},
	which, for $ x \in \HH $ outputs $ P_C ( x ) $, the orthogonal projection of $ x $ onto $ C $.
	\\
	\\
	We expect to derive high complexity for algorithm using the projection oracle as the only information available on $ C $.
	\\
	As we already observed, at iteration $k\geq 1$, Algorithm \ref{alg:main} -- see Eqs. \eqref{eq:farias} and \eqref{eq:err.cri} --, demands the
	computation of a $\hat\sigma$-approximate solution of the monotone inclusion
	\eqref{eq:walt},
	which is equivalent to 
	\begin{align*}
		VIP(G_k,C)
	\end{align*}
	where $ G_k $ is as in \eqref{def:opg}.
	It is easy to check that $G_k(\cdot)$ is $L_k$-Lipschitz continuous and $\mu$-strongly monotone, where
	\begin{align} \lab{def:Lmu}
		L_k := \lambda_k\norm{F'(y_{k-1})} + 1\quad \mbox{and}\quad \mu = 1.
	\end{align}
	%
	%
	Summarizing, the subproblems of our main algorithm are variational inequalities for (affine) strongly-monotone and Lipschitz continuous operators. Hence, for solving \eqref{eq:walt} (assuming that the orthogonal projection $P_C$ onto $C$ is available) one can apply the Tseng's forward-backward method~\cite{mon.sva-hpe.siam10, Tseng2000Modified} with initial guess $y_{k-1}$:
	\begin{align} \lab{eq:tseng}
		\begin{cases}
			y^0 := y_{k-1},\quad  0 < s_k < 1/L_k\\[2mm]
			\mbox{For}\; j\geq 1:\\[2mm]
			\widetilde y^j = P_C\left(y^{j-1} - s_k\;G_k(y^{j-1})\right),\\[2mm]
			y^j = \widetilde y^j - s_k\left(G_k(\widetilde y^j) - G_k(y^{j-1})\right).
		\end{cases}
	\end{align}
	\bprop \lab{prop:cost}
	Consider the sequences $(y^j)$ and $(\widetilde y^j)$ as in \eqref{eq:tseng} and let
	\begin{align} \lab{eq:flu}
		\nu^j := \dfrac{1}{\lambda_k s_k}\left(y^{j-1} - \widetilde y^j - s_kG_k(y^{j-1})\right)\quad \mbox{for all}\quad j\geq 1,
	\end{align}
	where $G_k(\cdot)$ is as in \eqref{def:opg}. Then, for all $j\geq \widehat j_k$,
	\begin{align*}
		\begin{cases}
			\nu^j\in N_C(\widetilde y^j),\\[2mm]
			\left\|\lambda_k \left(F_{y_{k-1}}(\widetilde y^j) + \nu^j \right) + \widetilde y^j - x_{k-1}\right\|\leq
			\hat\sigma \norm{\widetilde y^j - y_{k-1}},
		\end{cases}
	\end{align*}
	where
	\begin{align} \lab{def:hatj}
		\widehat j_k :=1 + 
		\left\lceil \dfrac{2}{\omega_k}\max\left\{\log\left(2\min\left\{1/\sqrt{2s_k}, 1 + 1/\sqrt{1 - s_k^2L_k^2}\right\}\right),
		\log\left(\dfrac{2}{\hat\sigma s_k}\sqrt{\dfrac{1 + s_k L_k}{1 - s_k L_k}}\right)\right\}\right\rceil
	\end{align}
	and
	\begin{align} \lab{def:omega}
		\omega_k := \dfrac{2s_k(1 - s_k^2L_k^2)}{2s_k + 1 - s_k^2L_k^2}.
	\end{align}
	As a consequence, by choosing $s_k = 1/(2L_k)$, we conclude that the Tseng's forward-backward method \emph{(}as in \eqref{eq:tseng}\emph{)} finds a $\hat\sigma$-approximate solution of \eqref{eq:walt} in at most 
	\begin{align*}
		O\left(1 + \left\lceil (1 + L_k)\left(1 + \log^+ (L_k)\right)\right\rceil\right)
	\end{align*}
	iterations, where $L_k$ is as in \eqref{def:Lmu}.
	\eprop
	\bproo
	Let $d_{0, k} := \norm{y_{k-1} - y^*_k} > 0$, where $y^*_k\in C$ denotes the unique solution of \eqref{eq:walt}.
	From \cite[Proposition 6.1]{mon.sva-hpe.siam10} and \cite[Proposition 2.2]{Alves2016Regularized} we have
	\begin{align*}
		\norm{\widetilde y^j - y^*_k}\leq \min\left\{1/\sqrt{2s_k}, 1 + 1/\sqrt{1 - s_k^2L_k^2}\right\} 
		\left(1 - \omega_k\right)^{\frac{j - 1}{2}} \,d_{0,k}
		\quad \mbox{for all}\quad j\geq 1,
	\end{align*}
	where $\omega_k$ is as in \eqref{def:omega}.
	Hence, 
	\begin{align} \lab{eq:vasco}
		\norm{\widetilde y^j - y^*_k}\leq \dfrac{d_{0,k}}{2}\quad \mbox{for all}\quad j\geq j_k^*,
	\end{align}
	where
	\begin{align}\lab{def:js}
		j_k^* := 1 + \left\lceil \frac{2}{\omega_k}\log\left(2\min\left\{1/\sqrt{2s_k}, 1 + 1/\sqrt{1 - s_k^2L_k^2}\right\}\right) \right\rceil.
	\end{align}
	From \eqref{eq:vasco} and a simple argument based on the triangle inequality we obtain
	\begin{align} \lab{eq:cruz}
		\norm{\widetilde y^j - y_{k-1}}\geq \dfrac{d_{0,k}}{2}\quad \mbox{for all}\quad j\geq j_k^*.
	\end{align}
	Likewise, from \cite[Proposition 6.1]{mon.sva-hpe.siam10} and \cite[Proposition 2.2]{Alves2016Regularized} we find
	\begin{align*} 
		\left\|\lambda_k \left(F_{y_{k-1}}(\widetilde y^j) + \nu^j \right) + \widetilde y^j - x_{k-1}\right\|\leq \dfrac{1}{s_k}
		\sqrt{\dfrac{1 + s_k L_k}{1 - s_k L_k}} \left(1 - \omega_k\right)^{\frac{j - 1}{2}}\,d_{0,k}
		\quad \mbox{for all}\quad j\geq 1,
	\end{align*}
	where, for all $j\geq 1$, $\nu^j\in N_C(\widetilde y^j)$ is as in \eqref{eq:flu}.
	Hence,
	\begin{align} \lab{eq:ate}
		\left\|\lambda_k \left(F_{y_{k-1}}(\widetilde y^j) + \nu^j \right) + \widetilde y^j - x_{k-1}\right\|\leq 
		\hat\sigma \dfrac{d_{0,k}}{2}\quad \mbox{for all}\quad j\geq j_k^{**},
	\end{align}
	where
	\begin{align}\lab{def:jss}
		j_k^{**} := 1 + \left\lceil \dfrac{2}{\omega_k}\log\left(\dfrac{2}{\hat\sigma s_k}\sqrt{\dfrac{1 + s_k L_k}{1 - s_k L_k}}\right)\right\rceil.
	\end{align}
	Combining \eqref{def:js} -- \eqref{def:jss}, we get
	\begin{align}
		\left\|\lambda_k \left(F_{y_{k-1}}(\widetilde y^j) + \nu^j \right) + \widetilde y^j - x_{k-1}\right\|\leq
		\hat\sigma \norm{\widetilde y^j - y_{k-1}}\quad \mbox{for all}\quad j\geq \widehat j_k,
	\end{align}
	where $\widehat j_k$ is as in \eqref{def:hatj}. Note that the inclusion $\nu^j\in N_C(\widetilde y^j)$ follows easily from 
	\eqref{eq:tseng} and \eqref{eq:flu}. The last statement of the proposition follows by taking
	$s_k = 1/(2L_k)$ and using \eqref{def:hatj}.
	\eproo

\mgap

\section{Complexity analysis of Algorithm \ref{alg:main}}

In this section, we discuss the iteration-complexity of Algorithm \ref{alg:main}; 
the main results are Theorems \ref{th:main} and \ref{th:main02} below where pointwise 
$O\left(\frac{1}{\rho}\right)$
and ergodic $O\left(\frac{1}{\rho^{2/3}}\right)$ complexities are presented. Propositions \ref{prop:stepb}, \ref{prop:imp}
and \ref{prop:inc.theta} are technical results which will be needed. 
Proposition \ref{prop:special} shows that for the indexes to which the large-step condition 
$\lambda_k\norm{y_k - x_{k-1}}\geq \eta$ is satisfied, it follows that Algorithm \ref{alg:main} is a special instance
of Algorithm \ref{alg_gen} (see also Corollary \ref{cor:rate}). One the other hand, Proposition \ref{prop:iteb} deals with those
indexes for which $\lambda_k\norm{y_k - x_{k-1}}\geq \eta$ doesn't hold (see \eqref{eq:def.ab}--\eqref{eq:car.a} below).

Before stating the results of this section, we recall the notion of approximate solution for the problem \eqref{eq:probm} we consider in this paper: for a given tolerance $\rho>0$, find $y\in C$ such that either there exists $\nu\in N_C(y)$ such that
\begin{align}\lab{eq:app.sola}
\norm{F(y) + \nu}\leq \rho
\end{align}
or there exists $\varepsilon\geq 0$ and $v\in (F + N_C)^{\varepsilon}(y)$ such that
\begin{align}\lab{eq:app.solb}
\max \left\{\norm{v},\varepsilon\right\}\leq \rho.
\end{align}

To simplify the presentation, define for $k \geq 1$
\begin{align}\lab{eq:def.ek}
  e_k & := \lambda_k\left(F_{y_{k-1}}(y_k) + \nu_k\right) + y_k - x_{k-1}\;.
\end{align}

\mgap

\begin{proposition}\lab{prop:stepb}
Consider the sequences evolved by \emph{Algorithm \ref{alg:main}}. Then, for all $k\geq 1$,
\begin{align*}
\norm{y_k - y_{k-1}}\leq \dfrac{1}{1-\hat\sigma}\norm{\lambda_k \left(F(y_{k-1})+\nu_{k-1}\right) + y_{k-1} - x_{k-1}}.
\end{align*}
\end{proposition}
\bproo
First note that if $y_k$ is computed as in line \ref{eq:n.cond} of Algorithm \ref{alg:main}, then the inequality holds trivially. 
Assume now that $y_k$ is computed as in line \ref{eq:do.newton} of Algorithm \ref{alg:main}, so that \eqref{eq:err.cri} holds.
In this case, note that from Algorithm \ref{alg:main}'s definition we have that $\nu_k\in N_C(y_k)$ for all $k\geq 0$. Hence, in view of the monotonicity of $N_C$, we obtain
\begin{align}\lab{eq:lego}
\inner{\nu_k - \nu_{k-1}}{y_k - y_{k-1}}\geq 0\qquad \forall k\geq 1.
\end{align}
Now by multiplying the inequality in \eqref{eq:lego} by $\lambda_k>0$ 
and using \eqref{eq:def.ek}
we have
\begin{align}\lab{eq:mouse}
\inner{e_k + x_{k-1} - y_k - \lambda_k F_{y_{k-1}}(y_k) - \lambda_k\nu_{k-1}}{y_k - y_{k-1}}\geq 0.
\end{align}

Note also that using \eqref{eq:def.ek} and the inequality in \eqref{eq:err.cri} we also find
\begin{align}\lab{eq:lego3}
\norm{e_k}\leq \hat\sigma \norm{y_k - y_{k-1}}
\qquad \forall k\geq 1.
\end{align}
After some simple algebraic manipulations with \eqref{eq:mouse} we obtain
\begin{align}\lab{eq:lego2}
\nonumber
\normq{y_k - y_{k-1}}&\leq \inner{e_k + x_{k-1} - y_{k-1}- \lambda_k F_{y_{k-1}}(y_k) - \lambda_k\nu_{k-1}}{y_k - y_{k-1}}\\[2mm]
\nonumber
&=\inner{e_k + x_{k-1} - y_{k-1}- \lambda_k (F(y_{k-1})+\nu_{k-1})}{y_k - y_{k-1}}\\[2mm]
\nonumber
&\hspace{2cm}- \lambda_k \inner{F'(y_{k-1})(y_k - y_{k-1})}{y_k - y_{k-1}}
\\[2mm]
\nonumber
&\leq \inner{e_k + x_{k-1} - y_{k-1}- \lambda_k (F(y_{k-1})+\nu_{k-1})}{y_k - y_{k-1}}\\[2mm]
&\leq \big(\norm{e_k} + \norm{x_{k-1} - y_{k-1}- \lambda_k (F(y_{k-1})+\nu_{k-1})}\big)\norm{y_k - y_{k-1}},
\end{align}
where we also used the definition of $F_{y_{k-1}}(y_k)$ -- see \eqref{eq:def.lin} -- and the fact that $F'(y_{k-1})\geq 0$ (see, e.g., \cite[Proposition 12.3]{roc.wet-book}).

To finish the proof of the proposition, note that from \eqref{eq:lego3} and \eqref{eq:lego2} we find
\begin{align*}
\norm{y_k - y_{k-1}}\leq \hat\sigma\norm{y_k - y_{k-1}} + \norm{x_{k-1} - y_{k-1}- \lambda_k (F(y_{k-1})+\nu_{k-1})},
\end{align*}
which is clearly equivalent to the desired inequality.
\eproo

\mgap

\bprop \lab{prop:imp}
Consider the sequences evolved by \emph{Algorithm \ref{alg:main}}, let $\hat\theta>0$ be as in \eqref{eq:def.hth} and assume that at iteration $k\ge1$, $y_k$ and $\nu_k$ are computed as in line \ref{eq:do.newton} of \emph{Algorithm \ref{alg:main}}, so that  \eqref{eq:err.cri} holds. Then:
\begin{align*}
\dfrac{\lambda_k L}{2}\norm{\lambda_k \left(F(y_{k-1}) + \nu_{k-1}\right) + y_{k-1} - x_{k-1}}\leq 
\theta\;\;\Rightarrow\;\;
\dfrac{\lambda_k L}{2}\norm{\lambda_k \left(F(y_k) + \nu_k\right) + y_k - x_{k-1}}\leq 
\hat\theta.
\end{align*}
\eprop
\bproo
Assume that $y_k$ is computed as in line \ref{eq:do.newton} of Algorithm \ref{alg:main} and
\begin{align}\lab{eq:coffee}
\frac{\lambda_k L}{2}\norm{\lambda_k \left(F(y_{k-1}) + \nu_{k-1}\right) + y_{k-1} - x_{k-1}}\leq 
\theta.
\end{align}  
Letting $e_k$ be as in \eqref{eq:def.ek}
and using its definition, some simple algebraic manipulations 
and the inequality in \eqref{eq:err.cri}
 we find
\begin{align*}
\norm{\lambda_k \left(F(y_k) + \nu_k\right) + y_k - x_{k-1}} &= 
               \norm{e_k + \lambda_k\left(F(y_k) - F_{y_{k-1}}(y_k)\right)}\\[2mm]
               &\leq \norm{e_k} + \lambda_k\norm{F(y_k) - F_{y_{k-1}}(y_k)}\\[2mm]
               &\leq \hat\sigma \norm{y_k - y_{k-1}} + \dfrac{\lambda_k L}{2}\normq{y_k - y_{k-1}},
\end{align*}
where in the second inequality we also used Lemma \ref{lm:lip} below. 

As a consequence, from Proposition \ref{prop:stepb},
\begin{align*}
\norm{\lambda_k \left(F(y_k) + \nu_k\right) + y_k - x_{k-1}} &\leq 
\dfrac{\hat\sigma}{1-\hat\sigma}\norm{\lambda_k \left(F(y_{k-1}) + \nu_{k-1}\right) + y_{k-1} - x_{k-1}}\\[2mm]
& \hspace{2cm}+\dfrac{1}{(1-\hat\sigma)^2}\dfrac{\lambda_k L}{2}
\norm{\lambda_k \left(F(y_{k-1}) + \nu_{k-1}\right) + y_{k-1} - x_{k-1}}^2.
\end{align*}
By multiplying the latter inequality by $\frac{\lambda_k L}{2}$, using the assumption 
\eqref{eq:coffee} and \eqref{eq:def.hth} we obtain
\begin{align*}
\dfrac{\lambda_k L}{2}\norm{\lambda_k F(y_k) + y_k - x_{k-1}} &\leq \theta\left(\dfrac{\hat\sigma}{1-\hat\sigma}
+ \dfrac{\theta}{(1-\hat\sigma)^2}\right)\\[2mm]
 & = \hat \theta,
\end{align*}
which finishes the proof of the proposition.
\eproo

\mgap

\bprop \lab{prop:inc.theta}
Consider the sequences evolved by \emph{Algorithm \ref{alg:main}} and let $\hat\theta$ be as in \eqref{eq:def.hth}. The following holds:
\begin{itemize}
\item[\emph{(a)}] For all $k\geq 1$,
\begin{align}\lab{eq:inc.theta}
\dfrac{\lambda_k L}{2}\norm{\lambda_k \left(F(y_{k-1}) + \nu_{k-1}\right) + y_{k-1} - x_{k-1}}\leq 
\theta.
\end{align}
\item[\emph{(b)}] For all $k\geq 1$,
\begin{align}
\dfrac{\lambda_k L}{2}\norm{\lambda_k \left(F(y_k) + \nu_k\right) + y_k - x_{k-1}}\leq 
\hat\theta.
\end{align}
\end{itemize}
\eprop
\bproo 
(a) Let us proceed by induction on $k\geq 1$. Note first that \eqref{eq:inc.theta} is true for $k=1$ in view of 
the facts that $y_0=x_0$, $\nu_0 =0$ and $\frac{\lambda_1 L}{2}\norm{\lambda_1 F(y_0)}\leq \theta$ (see the
input in Algorithm \ref{alg:main}).
Assume now that \eqref{eq:inc.theta} holds for some $k\geq 1$. 
We claim that 
\begin{align}\lab{eq:sophia}
\dfrac{\lambda_k L}{2}\norm{\lambda_k \left(F(y_k) + \nu_k\right) + y_k - x_{k-1}}\leq \hat\theta.
\end{align}
Indeed, if $y_k$ and $\nu_k$ are given as in lines \ref{eq:n.cond} and \ref{eq:n.cond02} of Algorithm \ref{alg:main}, then we 
readily obtain \eqref{eq:sophia}. Otherwise, \eqref{eq:sophia} follows from the induction hypothesis and Proposition \ref{prop:imp}.

Next we will consider two cases: $\lambda_k\norm{y_k - x_{k-1}}\geq \eta$ and
$\lambda_k\norm{y_k - x_{k-1}} < \eta$. In the first case, using lines \ref{eq:ext.step} and \ref{eq:red.lam} of Algorithm \ref{alg:main} and \eqref{eq:sophia}, we find

\begin{align}\lab{eq:sophia2}
\nonumber
\dfrac{\lambda_{k+1} L}{2}\norm{\lambda_{k+1} \left(F(y_k) + \nu_k\right) + y_k - x_k} 
\nonumber
& = \dfrac{\lambda_{k+1} L}{2}\norm{\lambda_{k+1} \left(F(y_k) + \nu_k\right) + y_k - x_{k-1} + x_{k-1} - x_k}\\[2mm]
\nonumber
& =\dfrac{\lambda_{k+1} L}{2}\norm{\lambda_{k+1} \left(F(y_k) + \nu_k\right) + y_k - x_{k-1} + \tau\lambda_k
   \left(F(y_k) + \nu_k\right)}\\[2mm]
\nonumber
& = \dfrac{\lambda_{k+1} L}{2}\norm{\left(\lambda_{k+1} + \tau\lambda_k\right) \left(F(y_k) + \nu_k\right) 
+ y_k - x_{k-1}}\\[2mm]
\nonumber
& = (1-\tau)\dfrac{\lambda_k L}{2}\norm{\lambda_k \left(F(y_k) + \nu_k\right) + y_k - x_{k-1}}\\[2mm]
\nonumber
& \leq (1 - \tau)\hat \theta\\[2mm]
&\leq \theta,
\end{align}
where the latter inequality comes from the fact that $0<\tau<1$ and 
$\hat\theta < \theta$ (because of the assumption that $0<\theta<(1-\hat\sigma)(1 - 2\hat\sigma)$). 
Note that \eqref{eq:sophia2} finishes the induction argument in the case $\lambda_k\norm{y_k - x_{k-1}}\geq \eta$.

Let us now consider the case $\lambda_k\norm{y_k - x_{k-1}} < \eta$. Note first that, using the latter inequality, lines \ref{eq:null.step} and \ref{eq:inc.lam} of Algorithm \ref{alg:main} and \eqref{eq:sophia} combined with the induction hypothesis, we find
\begin{align*}
&\dfrac{\lambda_{k+1} L}{2}\norm{\lambda_{k+1} \left(F(y_k) + \nu_k\right) + y_k - x_k}\\[3mm] 
 = &\dfrac{\lambda_{k+1} L}{2}\norm{\lambda_{k+1} \left(F(y_k) + \nu_k\right) + y_k - x_{k-1}}\\[2mm]
 = &\dfrac{\lambda_{k+1} L}{2}\Big\|\dfrac{\lambda_{k+1}}{\lambda_k}\big[\lambda_k \left(F(y_k) + \nu_k\right) + y_k - x_{k-1}\big] \\[2mm]
&\hspace{2cm} + \left(1 - \dfrac{\lambda_{k+1}}{\lambda_k}\right)(y_k - x_{k-1})\Big\|\\[2mm]
\leq &\left(\dfrac{\lambda_{k+1}}{\lambda_k}\right)^2\dfrac{\lambda_k L}{2}
 \norm{\lambda_k \left(F(y_k) + \nu_k\right) + y_k - x_{k-1}}\\[2mm] 
&\hspace{2cm}+ \dfrac{\lambda_{k+1} L}{2}\left(\dfrac{\lambda_{k+1}}{\lambda_k} - 1\right)\norm{y_k - x_{k-1}}\\[2mm]
< &\left(\dfrac{\lambda_{k+1}}{\lambda_k}\right)^2 \hat \theta + \dfrac{\lambda_{k+1} L}{2}
\left(\dfrac{\lambda_{k+1}}{\lambda_k} - 1\right)
\dfrac{\eta}{\lambda_k}\\[2mm]
= &\left(\dfrac{\lambda_{k+1}}{\lambda_k}\right)^2 \hat \theta + \dfrac{\lambda_{k+1}}{\lambda_k}
\left(\dfrac{\lambda_{k+1}}{\lambda_k} - 1\right)
\dfrac{\eta L}{2}\\[2mm]
 = &\dfrac{\hat\theta}{(1 - \tau)^2} + \dfrac{\tau}{(1-\tau)^2}\dfrac{\eta L}{2}.
\end{align*}
As a consequence, to complete the induction argument, it is sufficient to show that
\begin{align*}
\dfrac{\hat\theta}{(1 - \tau)^2} + \dfrac{\tau}{(1-\tau)^2}\dfrac{\eta L}{2} = \theta,
\end{align*}
or, equivalently,
\begin{align*}
(1-\tau)^2\theta - \dfrac{\eta L}{2}\tau - \hat \theta = 0.
\end{align*}
Note that the latter identity is clearly equivalent to
\begin{align*}
\theta \tau^2 - \left(2\theta + \dfrac{\eta L}{2}\right)\tau + \theta -\hat\theta = 0.
\end{align*}
Hence, the result follows from \eqref{eq:def.q} and \eqref{eq:root.q}. This completes the induction argument and finishes the proof of item (a).

\mgap

(b) This is a direct consequence of item (a), Proposition \ref{prop:imp} and lines \ref{eq:n.cond} and \ref{eq:n.cond02}
of Algorithm \ref{alg:main}.
\eproo

\mgap
\mgap

For every $j\geq 1$, define
\begin{align}
  \lab{eq:def.ab}
  \begin{aligned}
    A_j & =\left\{1\leq k \leq j\;|\;\; \text{$\lambda_k\norm{y_k - x_{k-1}}\geq \eta$\; holds at iteration $k$}\right\},
	  & a_j =\# A_j\,,\\
    B_j & =\{1\leq k\leq j\;|\;\; 
   \text{$\lambda_k\norm{y_k - x_{k-1}} < \eta$\; holds at iteration $k$}\},
    & b_j =\# B_j\,,
  \end{aligned}
\end{align}
\begin{align}
 \lab{eq:uni.ab}
 A =\bigcup_{j=1}^\infty A_j\;\;\mbox{and}\;\; B =\bigcup_{j=1}^\infty B_j.
\end{align}
We also define
\begin{align}
\lab{eq:s70}
  K=\set{k\geq 1\;|\; k\leq \# A},\quad
  i_0=0,\quad i_k=\text{$k$-th element of }A,
\end{align}
and note that 
\begin{align}
 \lab{eq:car.a}
  i_0<i_1<i_2<\cdots,\qquad A=\set{i_k\;|\; k\in K}.
\end{align}

\mgap

\brema \lab{rm:rm}
\emph{
We claim that
if $(x_k)$ is generated by Algorithm \ref{alg:main} then
\begin{align}
\lab{eq:ixk}
 x_{i_k-1}=x_{i_{k-1}}\qquad \forall k\in K,
\end{align}
where $K$ is as in \eqref{eq:s70} -- see also \eqref{eq:car.a}.
Indeed, for any $k\in K$, by \eqref{eq:uni.ab} and \eqref{eq:car.a} we have 
$\{k\geq 1\,|\, i_{k-1}< k <i_k\}\subset B$. Consequently, by the definition 
of $B$ as in \eqref{eq:uni.ab},  and line \ref{eq:null.step} we conclude that 
$x_i=x_{i_{k-1}}$ whenever $i_{k-1}\leq i<  i_k$. 
As a consequence, we obtain that \eqref{eq:ixk} follows
from the fact that $i_{k-1}\leq i_k-1<i_k$.}
\erema

\mgap

\bprop \lab{prop:special}
Consider the sequences evolved by \emph{Algorithm \ref{alg:main}} and let the set $K$ be as in \eqref{eq:s70} 
and \eqref{eq:car.a}.
Then, for all $k\in K$,
\begin{align}
  %
  \lab{eq:stick00}
  &\nu_{i_k}\in N_C(y_{i_k}),\\[2mm]
  \lab{eq:stick01}
   &\left\|\lambda_{i_k} \left(F(y_{i_k}) + \nu_{i_k}\right) + y_{i_k} - x_{i_{k-1}}\right\|\leq \sigma\norm{y_{i_k} - x_{i_{k-1}}},\\[2mm]
   \lab{eq:stick02}
   &\lambda_{i_k}\norm{y_{i_k} - x_{i_{k-1}}}\geq \eta,\\[2mm]
   \lab{eq:stick03}
   & x_{i_k} = x_{i_{k-1}} - \tau\lambda_{i_k} \left(F(y_{i_k}) + \nu_{i_k}\right),
  %
  \end{align}
  where
 \begin{align}\lab{eq:def.sigma}
 \sigma:=\dfrac{2\hat\theta}{\eta L}<1.
 \end{align}
 As a consequence, the sequences $(x_{i_k})_{k\in K}$, $(y_{i_k})_{k\in K}$ and $(\lambda_{i_k})_{k\in K}$ are generated
 by \emph{Algorithm \ref{alg:rHPE}} \emph{(}with $N=\# A$\emph{)} for solving \eqref{eq:probm}.
\eprop
\bproo
Note first that \eqref{eq:stick00} is a direct consequence of Algorithm \ref{alg:main}'s definition and the definition 
of $A$ as in \eqref{eq:uni.ab}.
Moreover, note that \eqref{eq:stick02} and \eqref{eq:stick03} follow directly from \eqref{eq:def.ab}--\eqref{eq:car.a}
and the assumption that $k\in K$ (see also line \ref{eq:ext.step} in Algorithm \ref{alg:main}).
Let's now prove the inequality \eqref{eq:stick01}.
From Proposition \ref{prop:inc.theta}(b) we obtain
\[
\dfrac{\lambda_{i_k}L}{2}\left\|\lambda_{i_k} \left(F(y_{i_k}) + \nu_{i_k}\right) + y_{i_k} - x_{i_k-1}\right\|\leq \hat\theta,
\]
which when combined with \eqref{eq:ixk}, \eqref{eq:def.sigma} and \eqref{eq:stick02} yields
\begin{align*}
\left\|\lambda_{i_k} \left(F(y_{i_k}) + \nu_{i_k}\right) + y_{i_k} - x_{i_{k-1}}\right\|&\leq \dfrac{2\hat\theta}{\lambda_{i_k} L}\\[2mm]
    & = \sigma\dfrac{\eta}{\lambda_{i_k}}\\[2mm]
    & \leq \sigma \norm{y_{i_k} - x_{i_{k-1}}},
\end{align*}
which gives exactly \eqref{eq:stick01}. The last statement in the proposition is a direct consequence of 
\eqref{eq:stick01}--\eqref{eq:stick03}, Algorithm \ref{alg:rHPE}'s definition and the definition of $K$ as in \eqref{eq:s70}.
\eproo

\mgap

\bcoro \lab{cor:rate}
Consider the sequences evolved by \emph{Algorithm \ref{alg:main}}, and the ergodic sequences as in \eqref{eq:def.erg1}--\eqref{eq:def.erg3} -- for the sequences $(y_{i_k})_{k\in K}$ and $(\lambda_{i_k})_{k\in K}$, where the set $K$ is as in \eqref{eq:s70}. 
Let $d_0$ denote the distance of $x_0$ to the solution set $(F + N_C)^{-1}(0)\neq \emptyset$ of \eqref{eq:probm}. 
Let also $0<\tau<1$ and $0\leq\sigma<1$ be as in \eqref{eq:def.tau.alg} and \eqref{eq:def.sigma}, respectively.

The following hold for all $k\in K$:
\begin{itemize}
\item[\emph{(a)}] There exists $1\leq j\leq k$ such that $\nu_{i_j}\in N_C(y_{i_j})$ and
\begin{align*}
\norm{F(y_{i_j}) + \nu_{i_j}}\leq \dfrac{d_0^2}{\tau \eta(1-\sigma)\,k},
\end{align*}
\item[\emph{(b)}] $v_k^a\in (F + N_C)^{\varepsilon_k^a}(y_k^a)$ and
\begin{align*}
\norm{v_k^a}&\leq \dfrac{2d_0^2}{\tau^{3/2}\eta\sqrt{1-\sigma^2}\,k^{3/2}},\\[3mm]
\varepsilon_k^a&\leq \dfrac{2d_0^3}{\tau^{3/2}\eta(1-\sigma^2)\,k^{3/2}}.
\end{align*}
\end{itemize}
\ecoro
\bproo
The proof is a direct consequence of the last statement in Proposition \ref{prop:special} and Theorem \ref{th:rhpe}.
\eproo

\mgap

\bprop \lab{prop:iteb}
Consider the sequences evolved by \emph{Algorithm \ref{alg:main}} and let the set $B$ be as in \eqref{eq:def.ab} and
\eqref{eq:uni.ab}. Then, for all $k\in B$,
\begin{align}\lab{eq:res.b}
\nu_k\in N_C(y_k)\;\;\mbox{and}\;\;\norm{F(y_k) + \nu_k}\leq \dfrac{2\hat\theta/L + \eta}{\lambda_k^2}, 
\end{align}
where $\eta > 0$ is as in the input of \emph{Algorithm \ref{alg:main}} and $\hat\theta>0$ is as in \eqref{eq:def.hth}. 
\eprop
\bproo
Note first that the inclusion in \eqref{eq:res.b} follows directly from Algorithm \ref{alg:main}'s definition.
Now, from Proposition \ref{prop:inc.theta}(b), the definition of $B$ and the triangle inequality, we obtain
\begin{align*}
\norm{\lambda_k \left(F(y_k) + \nu_k\right)} &\leq \norm{\lambda_k \left(F(y_k) + \nu_k\right) + y_k - x_{k-1}}+\norm{y_k - x_{k-1}}\\[2mm]
   &\leq \dfrac{2\hat\theta}{\lambda_k L} + \dfrac{\eta}{\lambda_k},
\end{align*}
which clearly gives the desired inequality in \eqref{eq:res.b}.
\eproo

\mgap

\begin{lemma} \lab{lm:geo}
Let $(\lambda_k)$ be generated by \emph{Algorithm \ref{alg:main}} and let $(a_k)$ and $(b_k)$ be as in \eqref{eq:def.ab}.
Then, for all $k\geq 1$,
\begin{align}
\lambda_{k+1} = (1 - \tau)^{a_{k} - b_{k}}\lambda_1.
\end{align}
\end{lemma}
\bproo
The result follows from lines \ref{eq:red.lam} and \ref{eq:inc.lam} of Algorithm \ref{alg:main} and the definitions of
$a_k$ and $b_k$ as in \eqref{eq:def.ab}.
\eproo

\mgap

Next is our first result on the iteration-complexity of Algorithm \ref{alg:main}. As we mentioned before, it shows, as presented in items (a) and
(b), respectively, up to an additive logarithmic term, 
\emph{pointwise} $O\left(\frac{1}{\rho}\right)$
and \emph{ergodic} $O\left(\frac{1}{\rho^{2/3}}\right)$ iteration-complexities, where $\rho>0$ is a given tolerance.

\mgap

\begin{theorem}[{\bf First result on the iteration-complexity of Algorithm \ref{alg:main}}] \lab{th:main}
Let $\rho>0$ be a given tolerance and let $\hat\theta>0$, $\eta>0$, $0<\tau<1$ and $\lambda_1>0$ be as in
the input of \emph{Algorithm \ref{alg:main}}. 
Let also $d_0$ be the distance of $x_0$ to the solution set $(F + N_C)^{-1}(0)\neq \emptyset$ of \eqref{eq:probm}, let
$L>0$ be as in \eqref{eq:mon.lip}
and let $0<\sigma<1$ be as in \eqref{eq:def.sigma}.
The following holds:
\begin{itemize}
\item[\emph{(a)}] \emph{Algorithm \ref{alg:main}} finds $y\in \HH$ and $\nu\in N_C(y)$ such that
\begin{align}\lab{eq:pointw}
\norm{F(y) + \nu}\leq \rho
\end{align}
in at most 
\begin{align}\lab{eq:def.m}
 \left\lceil{\left(\dfrac{2}{\tau \eta(1-\sigma)}\right)\dfrac{d_0^2}{\rho}}\right\rceil + 
 \left\lceil{\dfrac{1}{2\tau}\log^+\left(\dfrac{\eta + 2\hat\theta/L}{\lambda_1^2\rho}\right)}\right\rceil
\end{align}
iterations.
\item[\emph{(b)}] \emph{Algorithm \ref{alg:main}} finds a triple $(y,v,\varepsilon)\in \HH\times \HH\times \R_+$  such
that
\begin{align}\lab{eq:erg}
v\in (F + N_C)^{\varepsilon}(y)\;\;\emph{and}\;\; \max\left\{\norm{v},\varepsilon\right\}\leq \rho
\end{align}
in at most
\begin{align}\lab{eq:def.mt}
\begin{aligned}
& \max\left\{\left\lceil{\left(\dfrac{2\sqrt[3]{4}}
      {\tau \eta^{2/3}(1-\sigma^2)^{1/3}}\right)\left(\dfrac{d_0^2}{\rho}\right)^{2/3}}\right\rceil, 
      \left\lceil{\left(\dfrac{2\sqrt[3]{4}}
      {\tau \eta^{2/3}(1-\sigma^2)^{2/3}}\right)\left(\dfrac{d_0^3}{\rho}\right)^{2/3}}\right\rceil\right\}\\[4mm]
&\hspace{3cm}+\left\lceil{\dfrac{1}{2\tau}\log^+\left(\dfrac{\eta + 2\hat\theta/L}{\lambda_1^2\,\rho}\right)}\right\rceil
\end{aligned}
\end{align}
iterations.
\end{itemize}
\end{theorem}
\bproo
(a) Define
\begin{align}\lab{eq:def.m1}
M_1 = \left\lceil{\left(\dfrac{1}{\tau \eta(1-\sigma)}\right)\dfrac{d_0^2}{\rho}}\right\rceil\;\;\mbox{and}\;\; 
M_2 = \left\lceil{\dfrac{1}{2\tau}\log^+\left(\dfrac{\eta + 2\hat\theta/L}{\lambda_1^2\,\rho}\right)}\right\rceil 
\end{align}
and note that, in view of \eqref{eq:def.m}, we have
\begin{align}\lab{eq:m1m2}
 M = 2M_1 + M_2,
\end{align}
where $M$ is defined as the quantity given in \eqref{eq:def.m}.
We will consider two cases: (i) $\#A \geq M_1$ and (ii) $\#A < M_1$, where $A$ is as in \eqref{eq:uni.ab}.
In the first case, it follows directly from Corollary \ref{cor:rate}(a) and the definition of $M_1$ as in \eqref{eq:def.m1} that
Algorithm \ref{alg:main} finds $y\in \HH$ satisfying \eqref{eq:pointw} in at most $M_1$ iterations. Since $M\geq M_1$, we have the desired result.

Consider now the case (ii), i.e., $\#A < M_1$. Let $\ell\geq 1$ be such that $\#A = a_\ell = a_k$ for all $k\geq \ell$ (see
\eqref{eq:def.ab} and \eqref{eq:uni.ab}). Hence, if $b_{k-1}\geq M_1 + M_2$, for some $k\geq \ell + 1$, then it follows that
\begin{align}\lab{eq:bma}
b_{k-1} - a_{k-1} = b_{k-1} - \#A > b_{k-1} - M_1 \geq M_2,
\end{align}
where in the latter inequality we also used 
the assumption that $b_{k-1}\geq M_1 + M_2$.

To conclude that proof of (a), note that \eqref{eq:bma}, Lemma \ref{lm:geo}, Proposition \ref{prop:iteb}, the definition of
$M_2$ as in \eqref{eq:m1m2} and some simple algebraic manipulation give that $\norm{F(y_k) + \nu_k}\leq \rho$.
Altogether, still in the case (ii), we proved that Algorithm \ref{alg:main} finds $y\in \HH$ satisfying \eqref{eq:pointw} in at most
$M_1 + (M_1 + M_2)$ iterations, that is, in at most $M$ iterations.

\mgap

(b) The proof follows the same outline of the proof in item (a), now using Corollary \ref{cor:rate}(b) instead of Corollary \ref{cor:rate}(a).
\eproo

\mgap

Next we present the results obtained in Theorem \ref{th:main} for the choices of $0<\theta<1$ and $\eta>0$ as given in \eqref{eq:def.te} below. Note that \eqref{eq:def.mL} and \eqref{eq:def.mtL} present, respectively, the pointwise and ergodic iteration-complexities of Algorithm \ref{alg:main} in terms of $0\leq \hat\sigma<1/2$, $L>0$, $\rho>0$ and 
$d_0$.

\mgap

\begin{theorem}[{\bf Second result on the iteration-complexity of Algorithm \ref{alg:main}}] \lab{th:main02}
Let $\rho>0$ be a given tolerance, let $0\leq \hat\sigma <1/2$, let
$L>0$ be as in \eqref{eq:mon.lip} and consider \emph{Algorithm 2} with the choices of
$0<\theta<1$ and $\eta>0$ as follows
\begin{align}\lab{eq:def.te}
\theta = \dfrac{(1-\hat\sigma)(1-2\hat\sigma)}{2}\;\;\mbox{and}\;\;\eta = \dfrac{4\theta}{L}.
\end{align}
\emph{(}the parameters $\hat\theta>0$, $0<\tau<1$ and $\lambda_1>0$ are considered as described in the input of \emph{Algorithm \ref{alg:main}}.\emph{)}
Let also $d_0$ be the distance of $x_0$ to the solution set $(F + N_C)^{-1}(0)\neq \emptyset$ of \eqref{eq:probm}.
The following holds:
\begin{itemize}
\item[\emph{(a)}] \emph{Algorithm \ref{alg:main}} finds $y\in \HH$ and $\nu\in N_C(y)$ satisfying \eqref{eq:pointw}
in at most 
\begin{align}\lab{eq:def.mL}
  \Gamma_{\mathrm p} := \left\lceil{\left(\dfrac{4}{1 - 2\hat\sigma}\right)^2\dfrac{L\,d_0^2}{\rho}}\right\rceil + 
 \left\lceil{(1 - \hat\sigma)\left(\dfrac{4}{1 - 2\hat\sigma}\right)\log^+\left(\dfrac{(1 - 2\hat\sigma)(5/2 - 2\hat\sigma)}{\lambda_1^2\,L\rho}\right)}\right\rceil
\end{align}
iterations.
\item[\emph{(b)}] \emph{Algorithm \ref{alg:main}} finds a triple $(y,v,\varepsilon)\in \HH\times \HH\times \R_+$
satisfying \eqref{eq:erg} in at most
\begin{align}\lab{eq:def.mtL}
\begin{aligned}
  \Gamma_{\mathrm e}& :=\max\left\{\left\lceil\left({\dfrac{16\sqrt[3]{4/3}\left(1 - \hat\sigma\right)^{1/3}}
      {(1 - 2\hat \sigma)^{5/3}}}\right)\left(\dfrac{L\,d_0^2}{\rho}\right)^{2/3}\right\rceil, 
      \left\lceil{\left(\dfrac{32\sqrt[3]{2/9}\left(1 - \hat\sigma\right)^{1/3}}
      {(1 - 2\hat \sigma)^{5/3}}\right)}\left(\dfrac{L\,d_0^3}{\rho}\right)^{2/3}\right\rceil\right\}\\[4mm]
&\hspace{3cm} +  \left\lceil{(1 - \hat\sigma)\left(\dfrac{4}{1 - 2\hat\sigma}\right)\log^+\left(\dfrac{(1 - 2\hat\sigma)(5/2 - 2\hat\sigma)}{\lambda_1^2\,L\rho}\right)}\right\rceil
\end{aligned}
\end{align}
iterations.
\end{itemize}
\end{theorem}
\bproo
(a) Direct use of the definition of $\theta$ as in \eqref{eq:def.te} and \eqref{eq:def.hth} yield
\begin{align}\lab{eq:tmth}
 \hat\theta = \frac{1 - 2\hat\sigma}{4}\;\;\mbox{and}\;\;\theta - \hat\theta = \frac{(1 - 2\hat\sigma)^2}{4},
\end{align}
which when combined with \eqref{eq:def.te} and \eqref{eq:tau.ineq} implies
\begin{align}\lab{eq:ineq.tsun}
\tau > \frac{\theta - \hat\theta}{4\theta} = \frac{1 - 2\hat\sigma}{8(1 - \hat\sigma)}.
\end{align}
Combining \eqref{eq:ineq.tsun} with 
the definition of $\eta$ as in \eqref{eq:def.te} we find
\begin{align}\lab{eq:ineq.te}
\tau \eta > \dfrac{(1 - 2\hat\sigma)^2}{4L}.
\end{align}
Using also \eqref{eq:def.sigma}, the definition of $\eta$ as in \eqref{eq:def.te}, the first identity in \eqref{eq:tmth} and
the assumption that $0\leq \hat\sigma < 1/2$ we also 
obtain $\sigma = \frac{1}{4(1 - \hat\sigma)}$ and so
\begin{align}\lab{eq:sigma.sun}
\sigma < \dfrac{1}{2}.
\end{align}
From \eqref{eq:ineq.te} and \eqref{eq:sigma.sun}, and some simple algebra, we obtain 
\begin{align}\lab{eq:dinamo}
\nonumber
\left(\dfrac{2}{\tau \eta(1-\sigma)}\right)\dfrac{d_0^2}{\rho}
&< \left(\dfrac{16 L}{(1 - 2\hat\sigma)^2}\right)\dfrac{d_0^2}{\rho}\\[2mm]
     & = \left(\dfrac{4}{1 - 2\hat\sigma}\right)^2\,\dfrac{L\,d_0^2}{\rho}.
\end{align}

On the other hand, direct use of \eqref{eq:def.te}, the first identity in \eqref{eq:tmth} and some simple calculations yield
\begin{align}\lab{eq:jazz}
\dfrac{\eta + 2\hat\theta/L}{\lambda_1^2\,\rho} = \dfrac{(1 - 2\hat\sigma)(5/2 - 2\hat\sigma)}{\lambda_1^2\,L\rho}.
\end{align}
Note now that the desired result follows directly from Theorem \ref{th:main}(a) combined with \eqref{eq:dinamo},
\eqref{eq:jazz} and \eqref{eq:ineq.tsun}.

\mgap

(b) Note first that using \eqref{eq:ineq.te}, \eqref{eq:def.te} and some simple computations we find
\begin{align}
 \tau \eta^{2/3} > \dfrac{(1 - 2\hat\sigma)^{5/3}}{\left(2^7(1 - \hat\sigma)\right)^{1/3}\,L^{2/3}},
\end{align}
which when combined with \eqref{eq:sigma.sun} yields
\begin{align}\lab{eq:val}
\dfrac{2\sqrt[3]{4}}{\tau \eta^{2/3}(1-\sigma^2)^{1/3}}\leq \dfrac{16\sqrt[3]{4/3}(1 - \hat\sigma)^{1/3}\,L^{2/3}}
{(1 - 2\hat\sigma)^{5/3}}.
\end{align}
By a similar reasoning we also have 
\begin{align}\lab{eq:val02}
\dfrac{2\sqrt[3]{4}}{\tau \eta^{2/3}(1-\sigma^2)^{2/3}}\leq \dfrac{32\sqrt[3]{2/9}(1 - \hat\sigma)^{1/3}\,L^{2/3}}
{(1 - 2\hat\sigma)^{5/3}}.
\end{align}

The desired result follows from \eqref{eq:val}, \eqref{eq:val02}, \eqref{eq:jazz}, \eqref{eq:ineq.tsun} and
Theorem \ref{th:main}(b).
\eproo

\mgap

\brema \lab{rem:hat}
\emph{Note that if we set $\hat\sigma = 0$ in Algorithm \ref{alg:main}, in which case the subproblem \eqref{eq:err.cri}
reduces to \eqref{eq:farias}, then the pointwise and ergodic iteration-complexities as in items (a) and (b) of Theorem \ref{th:main02}, respectively, are now given by
\begin{align}
 \left\lceil{\dfrac{16L\,d_0^2}{\rho}}\right\rceil + 
 \left\lceil{4\log^+\left(\dfrac{5}{2\lambda_1^2\,L\rho}\right)}\right\rceil
\end{align}
}
and
\begin{align}
\begin{aligned}
&\max\left\{\left\lceil 16\sqrt[3]{4/3}\left(\dfrac{L\,d_0^2}{\rho}\right)^{2/3}\right\rceil, 
      \left\lceil 32\sqrt[3]{2/9}\left(\dfrac{L\,d_0^3}{\rho}\right)^{2/3}\right\rceil\right\}\\[4mm]
&\hspace{3cm} +  \left\lceil{4\log^+\left(\dfrac{5}{2\lambda_1^2\,L\rho}\right)}\right\rceil.
\end{aligned}
\end{align}
\erema

\mgap
\mgap

A formal proof of next corollary is presented in Appendix~\ref{app:bigo}.

\begin{corollary}
  \lab{cr:bigo}
  Under the choice of the parameters as in \emph{Theorem~\ref{th:main02}},
  the pointwise and ergodic iteration-complexities of
  \emph{Algorithm~\ref{alg:main}} are, respectively,
  $O\left(\frac{1}{\rho}\right)$
  and
  $O\left(\frac{1}{\rho^{2/3}}\right)$.
\end{corollary}

\mgap

\newcommand{\textcolorpar}[2]{{\leavevmode\color{#1}#2}}
\section{Numerical experiments}
On our numerical experiments we focus on unconstrained convex-concave min-max problems. We note that any convex-concave min-max problem can be written as monotone variational inequality problem. For instance, we may rewrite $\min_x \max_y f(x, y)$ as a VIP by setting $F(x, y) = (\nabla f_x(x, y); -\nabla f_y(x, y))$. We consider the min-max problem
\begin{equation}\lab{eq:cubic_min_max_problem}
\min_{x\in \R^n} \max_{y\in \R^n} \frac{L}{6}\|x\|^3 + y^T(A x - b),
\end{equation}
where $A$ is an invertible $n\times n$ matrix and $b\in \R^n$. This is a convex-concave min-max problem, where the Hessian is $L$-Lipschitz. The unique saddle-point solution is given by $x_*=A^{-1} b$ and \linebreak $y_*=-\frac{L}{2} \|x_*\| (A^{-1})^T x_*$. 

We set $L=10^{-3}$ and generate $A$ as follows: we set $A = U S V^T$, where $U$ and $V$ are random orthogonal matrices, and $S$ is a diagonal matrix with diagonal entries $s_{ii} := 20^{-i/n}$. Our choice of parameters ensures that the condition number of $A$ is $20$. The entries of $b$ are generated independently from a standard normal distribution. 

\begin{figure}[t]
	\centering
	\includegraphics[width=\textwidth]{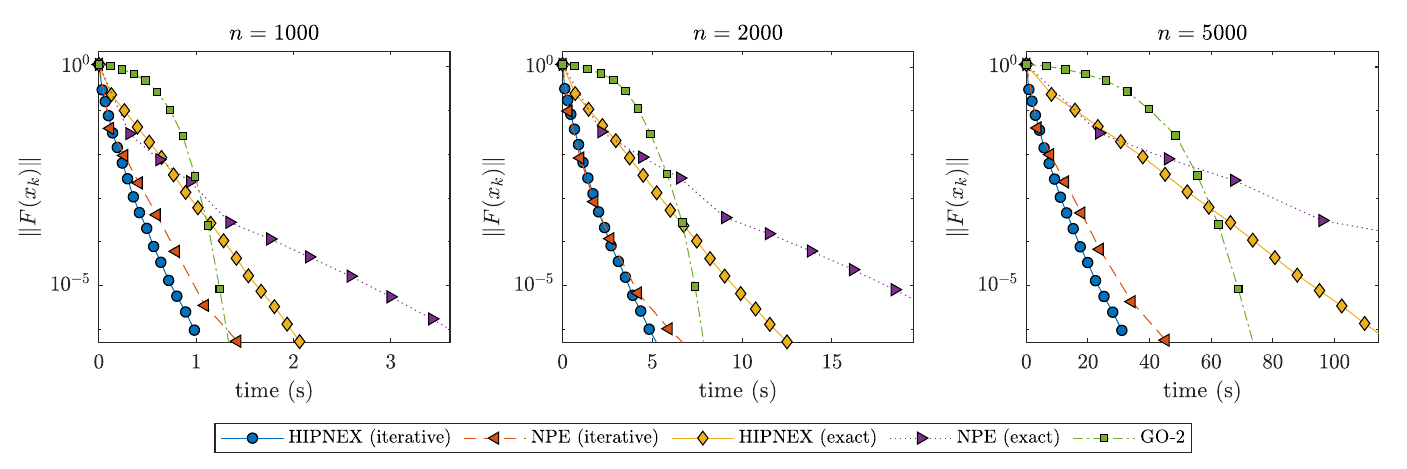}
	\caption{Time comparison between various optimization methods.}
	\lab{fig:numerical_results}
	\end{figure}

We compare the method we propose, HIPNEX\footnote{Acronym for Homotopy Inexact Proximal-Newton EXtragradient algorithm}, with the Newton Proximal Extragradient (NPE) \cite{mon.sva-newton.siam12} and the Second Order Generalized Optimistic (GO-2) \cite{jiang2022mokhtari} methods. The analysis for our method and NPE allow for solving sub-problems inexactly, as in Definition~\ref{def:app_sol}, therefore for each of these methods we compare two versions, one solving the linear system exactly (exact), using MATLAB default linear system solver, and the other using the MinRes \cite{paige1975solution} iterative solver (iterative), which is designed for solving symmetric indefinite linear systems. We compare these methods performance with the same initial point for $n=1000, 2000, 5000$. The entries of the initial points are generated independently from a standard normal distribution. We plot $\|F(x_k)\|$ for each iteration as a function of time in Figure~\ref{fig:numerical_results}, and provide further metrics on each method's performance in Table~\ref{table:numerical_results}. 
The numerical experiments were executed on a personal laptop with a AMD Ryzen 7 5800H CPU and 16.0GB of RAM. We share the code at \url{github.com/joaompereira/hipnex}.

%

\newcommand{\tento}[1]{\text{\footnotesize $\times10^{#1}$}}
\makeatletter
\newcommand*\ExpandableInput[1]{\input{#1} }
\makeatother
\begin{table}[t]
	\begin{footnotesize}
		\begin{center}
			\setlength\extrarowheight{2pt}			\setlength\tabcolsep{3pt}
			\begin{tabular}{|
					>{\centering\arraybackslash} m{16pt}|
					>{\centering\arraybackslash} m{88pt}|
					>{\centering\arraybackslash} m{34pt}|
					>{\centering\arraybackslash} m{38pt}|
					>{\centering\arraybackslash} m{45pt}|
					>{\centering\arraybackslash} m{36pt}|
					>{\centering\arraybackslash} m{36pt}|
					>{\centering\arraybackslash} m{36pt}|
					>{\centering\arraybackslash} m{42pt}|}
				\hline n & Method & time(s) & Iterations & $\|F(x_{\text{end}})\|$ & Linear Solves & $F$ evaluations & $J$ evaluations & Inner iterations\\
				\hline \multirow{5}*{$1000$} & HIPNEX (iterative) & $0.9934$ & $16$ & $9.59\tento{-07}$ & $16$ & $17$ & $16$ & $1870$ \\
\cline{2-9} & NPE (iterative) & $ 1.433$ & $7$ & $5.35\tento{-07}$ & $23$ & $14$ & $7$ & $2664$ \\
\cline{2-9} & HIPNEX (exact) & $ 2.068$ & $16$ & $5.22\tento{-07}$ & $16$ & $17$ & $16$ & \\
\cline{2-9} & NPE (exact) & $ 3.842$ & $10$ & $5.01\tento{-07}$ & $37$ & $20$ & $10$ & \\
\cline{2-9} & GO-2 & $ 1.379$ & $11$ & $1.53\tento{-07}$ & $11$ & $12$ & $11$ & \\
\hline \multirow{5}*{$2000$} & HIPNEX (iterative) & $ 5.315$ & $17$ & $4.40\tento{-07}$ & $17$ & $18$ & $17$ & $2010$ \\
\cline{2-9} & NPE (iterative) & $ 7.728$ & $7$ & $1.90\tento{-07}$ & $23$ & $14$ & $7$ & $2921$ \\
\cline{2-9} & HIPNEX (exact) & $ 12.51$ & $16$ & $5.17\tento{-07}$ & $16$ & $17$ & $16$ & \\
\cline{2-9} & NPE (exact) & $ 23.16$ & $10$ & $7.86\tento{-07}$ & $37$ & $20$ & $10$ & \\
\cline{2-9} & GO-2 & $ 8.047$ & $11$ & $1.76\tento{-07}$ & $11$ & $12$ & $11$ & \\
\hline \multirow{5}*{$5000$} & HIPNEX (iterative) & $ 31.04$ & $16$ & $9.48\tento{-07}$ & $16$ & $17$ & $16$ & $1853$ \\
\cline{2-9} & NPE (iterative) & $ 45.26$ & $7$ & $5.64\tento{-07}$ & $23$ & $14$ & $7$ & $2676$ \\
\cline{2-9} & HIPNEX (exact) & $   118$ & $16$ & $5.41\tento{-07}$ & $16$ & $17$ & $16$ & \\
\cline{2-9} & NPE (exact) & $ 252.7$ & $10$ & $5.98\tento{-07}$ & $37$ & $20$ & $10$ & \\
\cline{2-9} & GO-2 & $ 75.73$ & $11$ & $1.54\tento{-07}$ & $11$ & $12$ & $11$ & \\
\hline
			\end{tabular}
			
			\caption{Comparison between the algorithms in several metrics}
		\lab{table:numerical_results}
	\end{center} 
	
\end{footnotesize}
\end{table}

While comparing HIPNEX with NPE, we observe that the first is superior in several metrics, including the number of linear systems solved, the total number of inner iterations (for the iterative) versions, and consequently the computation time. Furthermore, we observe that the iterative methods were more advantageous for these test problems. Here we have not used pre-conditioners for speeding up the solution of iterative methods, however these should be considered in practical problems, where the Hessian is often badly-conditioned.

\appendix

\section{An auxiliary result}

A proof of the lemma below in finite-dimensional spaces can be found in \cite[Lemma 4.1.12]{den.sch-numbook.siam96}. 
A proof in general Hilbert spaces follows the same outline as in the finite-dimensional case.

\mgap

\begin{lemma}
 \lab{lm:lip}
 Let $F:\HH\to \HH$ be a continuously differentiable mapping
  with a $L$-Lipschitz continuous derivative $F'(\cdot)$, i.e., let $F$ be such that there exits $L\geq 0$ satisfying
\begin{align}
 \lab{eq:lip.app}
 \norm{F'(x) - F'(y)}\leq L\norm{x - y}\qquad  \forall x,y \in \HH.
\end{align}
Then,
  \[
   \norm{F(x) - F_y(x)}\leq \dfrac{L}{2}\norm{x - y}^2 \qquad \forall x,y\in \HH,
  \]
  where $F_y(x)$ is as in \eqref{eq:def.lin}.
\end{lemma}

\mgap

\section{Proof of Corollary~\ref{cr:bigo}}
 \lab{app:bigo}
  To estimate $\Gamma_{\mathrm p}$, first observe that for all $t>0$,
  \begin{align*}
 \lceil t \rceil  \leq 1 + t \quad \text{ and }\quad   \log^+(t) < t.
  \end{align*}
  Therefore,
\begin{align*}
  \left\lceil{\left(\dfrac{4}{1 - 2\hat\sigma}\right)^2\dfrac{L\,d_0^2}{\rho}}\right\rceil 
  \leq
  1+\left(\dfrac{4}{1 - 2\hat\sigma}\right)^2\dfrac{L\,d_0^2}{\rho}
\end{align*}
and
\begin{align*}
 \left\lceil{(1 - \hat\sigma)\left(\dfrac{4}{1 - 2\hat\sigma}\right)\log^+\left(\dfrac{(1 - 2\hat\sigma)(5/2 - 2\hat\sigma)}{\lambda_1^2\,L\rho}\right)}\right\rceil
 &\leq
1+ 4 
  \dfrac{(1 - \hat\sigma)(5/2 - 2\hat\sigma)}{\lambda_1^2\,L\rho}
  \\
 &\leq 1+ \dfrac{10}{\lambda_1^2\,L\rho}\,.
\end{align*}
Direct use of the above inequalities yields
\begin{align}
  \Gamma_{\mathrm p} \leq 2 + \dfrac{1}{\rho} 
  \left(
    \dfrac{16 L\,d_0^2}{(1 - 2\hat\sigma)^2}
    +
    \dfrac{10}{\lambda_1^2\,L}
  \right).
\end{align}
Hence, $\Gamma_{\mathrm p} = O(1/\rho)$.

To estimate $\Gamma_{\mathrm e}$, first observe that
\begin{align*}
  \max\{\lceil a \rceil , \lceil b \rceil\}
  \leq 
  \max\{ 1+a , 1+ b \}
  = 1 + \max\{ a ,  b \}.
\end{align*}
Hence,
\begin{align}
\nonumber
&\max\left\{\left\lceil
      \left({\dfrac{16\sqrt[3]{4/3}\left(1 - \hat\sigma\right)^{1/3}}
      {(1 - 2\hat \sigma)^{5/3}} }
      \right)
      \left(\dfrac{L\,d_0^2}{\rho}\right)^{2/3}\right\rceil, 
      \left\lceil{
            \left(\dfrac{32\sqrt[3]{2/9}\left(1 - \hat\sigma\right)^{1/3}}
      {(1 - 2\hat \sigma)^{5/3}}\right)
}\left(\dfrac{L\,d_0^3}{\rho}\right)^{2/3}\right\rceil\right\}\\[4mm]
\nonumber
&\hspace{3cm} \leq 
\\
\nonumber
& 1 + \max\left\{
  \left({\dfrac{16\sqrt[3]{4/3}\left(1 - \hat\sigma\right)^{1/3}}
      {(1 - 2\hat \sigma)^{5/3}}}\right)
    \left(\dfrac{L\,d_0^2}{\rho}\right)^{2/3}, 
          \left(\dfrac{32\sqrt[3]{2/9}\left(1 - \hat\sigma\right)^{1/3}}
      {(1 - 2\hat \sigma)^{5/3}}\right)
\left(\dfrac{L\,d_0^3}{\rho}\right)^{2/3}\right\}
\\[4mm]
\nonumber
&\hspace{3cm} =
\\
\lab{eq:tears}
& 1 +
    \left(
        \dfrac{
        32\left(1 - \hat\sigma\right)^{1/3}}
        {(1 - 2\hat \sigma)^{5/3}}
    \right)
      \left(\dfrac{L\,d_0^2}{\rho}\right)^{2/3} 
 \max\left\{
   \dfrac{1}{\sqrt[3]{6}}
   ,
   \sqrt[3]{ 2/9}
   \sqrt[3]{d_0^2}
\right\}\,.
%
\end{align}
Observe also that for any $t>0$ and $\alpha>0$,
\begin{align*}
 \log^+(t)= \dfrac{1}{\alpha}\log^+(t^\alpha)
 \leq \dfrac{t^\alpha}{\alpha}.
\end{align*}
Therefore,
\begin{align*}
  \log^+\left(
    \dfrac{
      (1 - 2\hat\sigma) (5/2 - 2\hat\sigma)}
    {\lambda_1^2\,L\rho}\right)
  & =\dfrac{3}{2} \log^+\left(
    \dfrac{(1 - 2\hat\sigma)(5/2 - 2\hat\sigma)}{\lambda_1^2\,L\rho}
 \right)^{2/3}
 \\
 & \leq \dfrac{3}{2} \left(\dfrac{(1 - 2\hat\sigma)(5/2 - 2\hat\sigma)}{\lambda_1^2\,L\rho}
 \right)^{2/3}
\end{align*}
and
\begin{align}
 \left\lceil{(1 - \hat\sigma)\left(\dfrac{4}{1 - 2\hat\sigma}\right)\log^+\left(\dfrac{(1 - 2\hat\sigma)(5/2 - 2\hat\sigma)}{\lambda_1^2\,L\rho}\right)}\right\rceil
\nonumber
 &\leq \\[2mm]
 \lab{eq:fears}
 & \hspace{-1.2cm}1+ (1 - \hat\sigma)\left(\dfrac{4}{1 - 2\hat\sigma}\right)
 \dfrac{3}{2}\left(\dfrac{(1 - 2\hat\sigma)(5/2 - 2\hat\sigma)}{\lambda_1^2\,L\rho}\right)^{2/3}.
\end{align}

Direct use of \eqref{eq:tears} and \eqref{eq:fears} then gives
\begin{align*}
\Gamma_e = O\left(\frac{1}{\rho^{2/3}}\right),
\end{align*}
which finishes the proof of Corollary \ref{cr:bigo}.




\def\cprime{$'$} \def\cprime{$'$}

\end{document}